\documentclass[a4paper,10pt]{article}
\usepackage[nointlimits]{amsmath}
\usepackage{amsfonts}
\usepackage{amssymb}
\usepackage{amsmath}
\usepackage{amsthm}
\usepackage{latexsym}
\usepackage[dvips]{graphicx}
\usepackage{layout}
\usepackage[english]{babel}
\usepackage{amscd}
\usepackage{color}
\newtheorem{lemma}     {Lemma}[section]
\newtheorem{thm}   [lemma]{Theorem}
\newtheorem{teorema1}   [lemma]{Theorem}
\newtheorem{prop}      [lemma]{Proposition}
\newtheorem{coro}    [lemma]{Corollary}
\newtheorem{cong1}      [lemma]{Conjecture}
\newtheorem{remark1}    [lemma]{Remark}

\numberwithin{equation}{section}

\renewcommand{\(}{\left(}        \renewcommand{\)}{\right)}
\renewcommand{\[}{\left[}        \renewcommand{\]}{\right]}

\newcommand{\sgn}{\mathrm{sign}}

\renewcommand{\i}{\infty}
\newcommand{\ra}{\rightarrow}

\renewcommand{\d}{\delta}

\renewcommand{\b}{\beta}
\newcommand{\eps}{\epsilon}

\newcommand{\ga}{\gamma}
\newcommand{\g}{\gamma}

\renewcommand{\o}{\omega}
\newcommand{\si}{\sigma}

\newcommand{\la}{\lambda}

\newcommand{\T}{{\mathcal T}}
\newcommand{\M}{{\mathcal M}}
\newcommand{\OO}{\mathcal O}
\newcommand{\F}{\mathcal F}
\newcommand{\p}{\mathcal P}
\renewcommand{\L}{\mathcal{L}}
\newcommand{\E}{\mathbb{E}}
\newcommand{\e}{\mathcal{E}}

\newcommand{\hc}{h_c}
\newcommand{\h}{h_N}
\newcommand{\ta}{\tau_{N,T}}

\title{Random hysteresis loops}

\author{{\normalsize{Gioia Carinci}}\\ \\\footnotesize{\emph{Dipartimento di Matematica Pura e Applicata, Universit\`{a} dell'Aquila}},\\
\footnotesize{\emph{via Vetoio 1, 67100 L'Aquila, Italy},}
\footnotesize{\emph{e-mail: gioia.carinci@univaq.it}}}
\date{}

\begin{document}
\maketitle

%\pagenumbering{arabic}
%\setcounter{equation}{0}

\begin{abstract}
Dynamical hysteresis  is a phenomenon which arises in
ferromagnetic systems below the critical temperature as a response
to adiabatic variations of the external magnetic field. We study
the problem in the context of the mean-field Ising model with
Glauber dynamics,  proving that for frequencies of the magnetic
field oscillations  of order $N^{-\frac 2 3}$, $N$ the size of the
system, the ``critical'' hysteresis loop becomes random.
\end{abstract}

\begin{abstract}
L'hyst\'{e}r\'{e}sis dynamique est un ph\'{e}nom\`{e}ne qu'on
observe dans les syst\`{e}mes ferromagn\'{e}tiques au-dessous de
la temperature critique, en r\'{e}ponse \`{a} des variations
adiabatiques du champ magn\'{e}tique ext\'{e}rieur. Nous
\'{e}tudions le probl\`{e}me dans le contexte du mod\'{e}le
d'Ising de  champ moyen avec la dynamique de Galuber, en montrant
que, pour des fr\'{e}quences d'oscillations du champ
magn\'{e}tique d'ordre de $N^{-2/3}$, avec $N$ la taille du
syst\`{e}me, la boucle d'hyst\'{e}r\'{e}sis ``critique'' devient
al\'{e}atoire.
\end{abstract}

%\pagenumbering{roman} \tableofcontents \newpage \pagenumbering{arabic}

\section{Introduction}

Hysteresis  appears when a time dependent magnetic field
$h=h(t)$ is applied to a ferromagnet  whose temperature is kept
fixed below the critical value.
 The origin of the phenomenon lies
in the fact that, at the equilibrium, at each value of the external magnetic field $h$
may not correspond a unique value of the magnetization $m$ of the
system. The value of $m(t)$ is, thus, not determined by $h(t)$
alone but also by the previous history of the input.

 The phenomenon has been widely studied and modelled. Most classical theories (see for example \cite{Be,BM,Vi})  consider 
hysteresis from a static point of view,  by modelling it through integral operators  not depending on the 
velocity of variation of the external input.\\
A dynamical approach to the study of the phenomenon  has been proposed for the first time by Rao et al. \cite{RKP} in the early nineties.
The new  theory aroused great interest and a number of experimental, numerical and theoretical works appeared on the argument in the last twenty years,
  investigating the response of the system to adiabatic oscillations of the magnetic field. 
 They analyse, in particular,   the dependence of shapes and areas of the hysteresis loops on  amplitude and frequency of the  input oscillations. 
   Most of these results are essentially numerical. Monte Carlo simulations have widely been used to study the hysteretic response 
  of a nearest-neighbor ferromagnetic  Ising model (see for instance \cite{AC,KNRW,KNRW1,NRS,NRS2,NRS3,RKP,ZDL}).
On the other hand, several   theoretical and numerical results are concerned with those  known as mean-field models (see \cite{AC,JGRM,RKP,TO}).
 In these models the dynamics is reduced to a single differential equation of the order parameter (the uniform magnetization $m(t)$). 
These equations govern  the dynamics of the magnetization in stochastic spin models in the limit of infinite system volume. Therefore they neglect both thermal fluctuations and finite system size effects.
A first rigorous analysis of  the effects of the stochastic fluctuations on the properties of the hysteresis cycles  has been carried out by B. Genz and N. Berglund in a series of papers of about ten years ago \cite{BG1,BG2,BG3}.
They model the thermal fluctuations by adding a stochastic noise to a mean-field type equation. They consider a Langevin equation with a Ginzburg-Landau potential:
        \begin{equation}
        \label{err1.1}
d x= (F(x)+h)dt + N^{-1/2} dw(t), \quad \quad \quad F(x)=x-x^3
     \end{equation}
where $w(t)$ is the standard brownian motion. 
We give to $N>0$
 the physical interpretation of the total number of spin sites
in a ferromagnetic system.
Then, in the large $N$ regime, equation \eqref{err1.1} can be thought of as a continuous counterpart 
of our Ising spin dynamics  (see Section \ref{section:Def}).
\\
In the present paper we shall study the
problem  for the Glauber process in the Curie-Weiss model, from which \eqref{err1.1}
is inspired. \\
 
\begin{figure}[htbp]
\centering
\includegraphics[width=100mm]{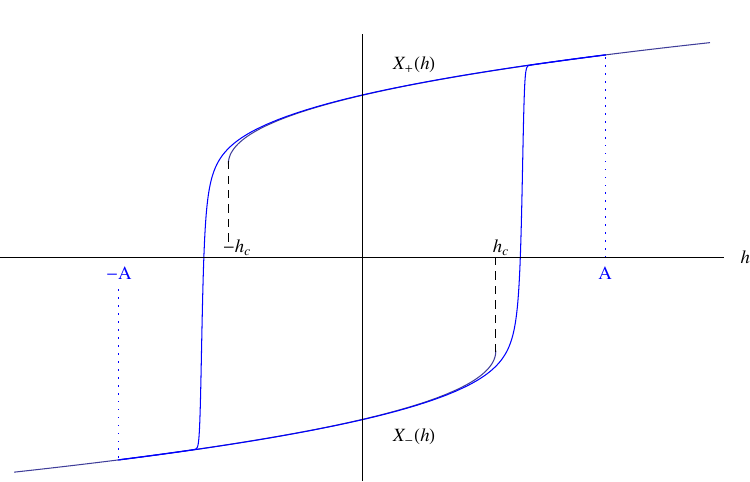}% "%" necessario
\caption{The picture shows the dependence of the two equilibrium branches $X_\pm(h)$ (black lines)
on the external magnetic field $h$. Adiabatic oscillations of the magnetic field of amplitude $A>h_c$ yield the typical hysteresis loop (blue line).}
\label{Hyste}
\end{figure}

 Let $h_c>0$ be the ``coercive magnetic field'' value, then for  $|h|\le
h_c$ the magnetization density of the ferromagnet may have two
equilibrium values, $X_+(h)$ and $X_-(h)$ (see Figure \ref{Hyste}). The upper
branch $X_+(h)$ continues past $h_c$ while it disappears for
$h<-h_c$; the opposite holds for the lower branch $X_-(h)$.  Let
us  apply, now,  a slowly oscillating magnetic field $h(t)$. We
denote, respectively, by $A$ and $\o$ the amplitude and the frequency
of the oscillations (we choose, for instance, $h(t)= -A \cos (\o
t)$). Let $m(t)$ be the magnetization observed at time $t$ and
choose initially $m(0)=X_+(h(0))$. In the adiabatic (quasi-static)
regime, where $\omega$ is very small, the following is observed.
If $A\le h_c$ then $m(t)\approx X_+(h(t))$ for any $t\ge 0$. If
$A> h_c$, $m(t)$ traces out the so called hysteresis loop, in the
sense that $m(t)\in\{ X_+(h(t)),X_-(h(t))\}$ (approximately),
jumping from the upper to the lower branch when $h(t)$ crosses $-
h_c$ and the opposite when $h(t)$ crosses $h_c$.
 A sharp statement (which avoids the above approximated
statements) can  be obtained in ``the adiabatic limit'' where
$\omega\to 0$.

 The pediod  of the magnetic field oscillations is of order $\o^{-1}$, thus, in
the adiabatic regime the natural time-scale of the dynamics is very long. 
In long time intervals
other phenomena may appear which in short time intervals are negligible and
which may invalidate the picture.   In the context of
\eqref{err1.1} $X_{\pm}(h)$ are identified with the locally stable
solutions of the stationary equation $F(x)=-h$.  If $h$ is
constant, say $h\in (0, h_c)$, then $X_-(h)$ is metastable and, on
a time interval which diverges exponentially with $N$ (as $N\to \infty$),
there is tunneling from $X_-(h)$ to $X_+(h)$. Thus, if $\omega$ is
exponentially small with $N$, the oscillations period is   exponentially
long with $N$, and then stochastic jumps between the two branches
occur, essentially perturbing the hysteresis loop. We intend to
consider a different regime for the frequency $\o$, i.e. we take
$\omega = N^{-\kappa}$, $\kappa>0$. We shall concentrate here on
the critical amplitude  case $A=h_c$. In such a case the deterministic
equation (i.e. \eqref{err1.1} without the brownian term) predicts
that the magnetization $m(t)$ tracks always the upper branch
$X_+(h(t))$, where it was initially. \cite{BG2} proves that, with
the addition of the stochastic effects,  there exists a critical
value for $\kappa$, $\kappa=\frac 23$. If $\kappa<\frac 23$ the
dynamics is still governed by the deterministic equation, i.e. the
magnetization tracks the upper branch, in the adiabatic limit.
Whereas, if $\kappa>\frac 23$ there is hysteresis, thus the
magnetization jumps to the lower branch as soon as $h=-h_c$ and
then back to the upper one when $h=h_c$ and so forth. 
In the present work we will prove  $\kappa=\frac 23$ to be the critical
value even in our Ising spin context. We shall
concentrate here on the critical case $\kappa=\frac 23$ which is
not covered by the analysis in \cite{BG2,BG3}.  
We will see that for $\kappa=\frac 23$ the hysteresis loop becomes
truly random.  There is a positive and not one probability
to leave the upper-lower branch at $\pm h_c$.
 Our  future aim is to extend our analysis to
the  Kac potential case by taking into account  spatial
effects.

\section{Definitions and results}
 \label{section:Def}

 {\bf The mean field Ising model.}  The configuration space is $\{-1,1\}^N$, $N\in \mathbb N$;
its elements are denoted by $\si=\{\si(i),i=1,..,N\}$, $\si(i)$
the spin at site $i$.  By
        \begin{equation}
        \label{err2.1}
m_N=m_N(\si):= \frac 1N \sum_{i=1}^N\si(i)
     \end{equation}
we denote the magnetization density of the configuration $\si$, so
that $m_N\in \mathcal{M}_N$,
          \begin{equation*}
\mathcal{M}_N:= \frac 1 N \Big\{ -N,-N+2,...,N -2 ,N\Big\}.
        \end{equation*}
The mean field  hamiltonian is
        \begin{equation*}
H_{h,N}(\si):=  N \Big( - \frac{m_N(\si)^2}{2}-hm_N(\si)\Big)
     \end{equation*}
and the mean field Gibbs measure at the inverse temperature
$\beta>0$ is the probability $G_{\beta,h,N}$ on  $\{-1,1\}^N$
given by
        \begin{equation*}
G_{\beta,h,N}(\si):= \frac {e^{-\beta
H_{h,N}(\si)}}{Z_{\beta,h,N}}
     \end{equation*}
where the partition function $Z_{\beta,h,N}$ is the normalization
factor.

\vskip.2cm
For an introduction to the mean field Ising model see Section 4.1 of \cite{Pr}.
\vskip.6cm

{\bf The Glauber dynamics.} A Glauber dynamics for the Ising
system is the Markov process on $\{-1,1\}^N$ with generator
        \begin{equation}
        \label{err2.5}
L f(\si):=\sum_{i=1}^N c(i, \si;h)\(  f(\si_i)-f(\si) \),
     \end{equation}
where $\si_i(j)=\si(j)$ for $i \neq j$ and $\si_i(i)=-\si(i)$;
$c(i, \si;h)>0$, the spin flip intensity at $i$, is given by the
formula
        \begin{equation*}
c(i, \si;h)= \frac{e^{-\beta[H_{h,N}(\si^{(i)})-
H_{h,N}(\si)]}}{e^{-\beta H_{h,N}(\si^{(i)})}+e^{-\beta
H_{h,N}(\si)}}
     \end{equation*}
with $\si^{(i)}$ the configuration obtained from $\si$ by flipping
the spin at $i$. For more details on the Glauber dynamics for mean field Ising systems 
see Section 5.1 of \cite{Pr}.
\vskip.2cm 

 $h=h(t)$ is a smooth function of time,
hence $\sigma(t)$ is a time non homogeneous Markov process.  Since the
hamiltonian depends on $\si$ via $m_N(\si)$, the process
$\{m_N(\si_t), t\ge 0\}$ is itself Markov with state space
$\mathcal M_N$ and generator $\L$ given by
        \begin{equation}\label{eq:L}
\L_h f(x):=c^+(x,h)\[ f(x+2 /N)- f(x)\]+c^-(x,h)\[ f(x- 2/ N)-
f(x)\]
\end{equation}
with
\begin{equation*}
 c^\pm(x,h)=\frac{N}2 (1\mp x)\, \hat c^\pm(x,h), \quad
\hat c^\pm(x,h)=\frac{e^{\pm \b[h+(x \pm 1/N)]}}{e^{-\b[h+(x\pm
1/N)]}+e^{\b[h+(x\pm 1/N)]}}\nonumber
    \end{equation*}

for $x \in \M_N$. When $h$ is time independent there is a unique
invariant measure (see Section 5.1.2 of \cite{Pr}) which is the marginal $\mu_{\beta,h,N}$ of
$G_{\beta,h,N}$ on the magnetization density $m_N$ defined in
\eqref{err2.1}. $\mu_{\beta,h,N}$ is then the probability on
$\mathcal M_N$ given by
        \begin{equation*}
\mu_{\beta,h,N}(x):= \frac {e^{-\beta N
\phi_{\beta,h,N}(x)}}{Z_{\beta,h,N}} \quad \quad x \in \M_N
     \end{equation*}
where
        \begin{equation*}
 \phi_{\beta,h,N}(x) := -\frac {x^2}{2} - hx -\frac{\mathcal S_{N}(x)}{\beta}
     \end{equation*}
and
       \begin{equation*}
e^{N \mathcal S_{N}(x)} := \text{\rm card}\Big(\si \in \{-1,1\}^N:
m_N(\si)=x\Big)
     \end{equation*}
If $x_N \in \M_N$, $x_N\to x\in[-1,1]$  as $N\to \infty$ then
$\phi_{\beta,h,N}(x_N) \to \phi_{\beta,h}(x)$ where
            \begin{equation*}
 \phi_{\beta,h}(x)= -\frac {x^2}{2} - hx -\frac{\mathcal S(x)}{\beta}
     \end{equation*}
and
                 \begin{equation*}
\mathcal S(x)=-\frac{1-x}2 \log \frac{1-x}2 - \frac{1+x}2 \log
\frac{1+x}2.
     \end{equation*}

\vskip.6cm

\begin{figure}[htbp]
\centering
\includegraphics[width=100mm]{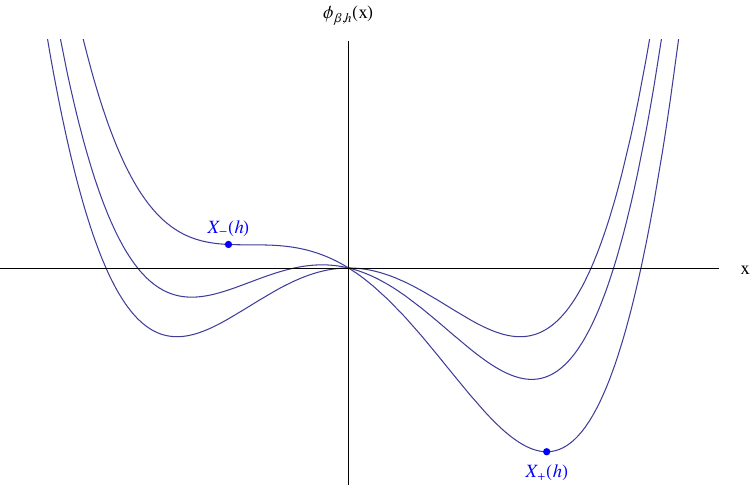}% "%" necessario
\caption{The picture shows some profiles of $\phi_{\beta,h}(x)$
according to different values of $h$, for $\b>1$.}\label{Phi}
\end{figure}

{\bf The mean field phase transitions.} For any $\beta\le 1$ and
any $h\in \mathbb R$ the mean field free energy density (see Section 4.1.2 of \cite{Pr})
$\phi_{\beta,h}(x)$ is a convex function of $x$ (absence of phase
transitions). If instead $\beta>1$ (see Figure \ref{Phi}) there is
$h_c>0$ such that, for any $|h|<h_c$, $\phi_{\beta,h}(x)$ is a
double well function of $x$ with local minima at $X_+(h)>X_-(h)$
and local maximum at $X_0(h)\in \big(X_-(h),X_+(h)\big)$;
$X_{\pm}(h)$ and $X_0(h)$ are solutions of the mean field
equation:
    \begin{equation*}
 x= \tanh\{\beta (x+h)\}
     \end{equation*}

\vskip.3cm

$X_+(h)$ is the absolute minimum for $h\ge 0$ and $X_-(h)$ for
$h\le 0$, then only at $h=0$ there are two absolute minima and
thus a phase transition; for $h\in (0,h_c)$, $X_+(h)$ is the only
pure phase while $X_-(h)$ is a metastable state, the opposite
holds for negative fields.  When $h\to -h_c$, $X_+(h)-X_0(h)\to 0$
and the limit $x_c:=X_+(-h_c)$ of $X_+(h)$ is an inflection point
for the function $\phi_{\beta,-h_c}(x)$.  By symmetry the
analogous picture describes $X_-(h)$ when $h\to h_c$.

\vskip.6cm

{\bf The macroscopic mean field dynamics.} The infinite volume
dynamics is governed by the ODE
\begin{equation}
    \label{eq:F}
\frac {dx}{dt}=F(x,h),\quad \quad F(x,h):= -x+\tanh\{ \b(x+h) \}
    \end{equation}
in the following sense. Let $m_N(t)$ be the process of generator
$\L_{h(t)}$ (see \eqref{eq:L}), $h(t)$ a smooth function of $t$,
which starts from $m^{0}_N\in \mathcal M_N$. We suppose that $m_N^0\to
x^0\in [-1,1]$ as $N\to \infty$ and denote by $\p_N$ the law of
$m_N(t), t\ge 0$. We have the following result.

\vskip.6cm

           \begin{thm}
           \label{thm:0}
With the above notation, for any $\d >0$ and any $T>0$,
     \begin{equation}
     \label{eq:th}
\lim_{N\to \infty}\p_N\bigg\{ \sup_{t \le T}\big|m_N(t)-x(t)\big|\ge \d \bigg\}=0
        \end{equation}
where $x(t)$ is the unique solution of
    \begin{equation}
    \label{eq:eqF}
\frac {dx}{dt}=F(x,h(t)),\quad x(0)=x^0
    \end{equation}
    \end{thm}

The proof of Theorem \ref{thm:0} is omitted. The proof in the case
of constant $h$ can be found, for instance, in Section  5.1.5
of \cite{Pr}, the proof  easily extends to the present case.

%\begin{figure}[!htbp]
%\centering
%\includegraphics[width=100mm]{m2}% "%" necessario
%\caption{The picture shows the stable branches $m_+(h(t))$ and
%$m_-(h(t))$ for $h(t)=-h_c \cos t$ connected by the unstable
%branch $m_0(h(t))$ (dashed line).}\label{m2}
%\end{figure}

\vskip.6cm

{\bf The adiabatic limit.} Let
 \begin{equation}
        \label{err2.18'}
h(t):= -h_c \cos t
     \end{equation}
we denote by $x_\omega(t)$ the solution of \eqref{eq:eqF}
with $h=h(\o t)$ and initial condition $x_\omega(0)=X_+(-h_c)$. We omit the proof that

\vskip.6cm

           \begin{thm}
           \label{thm:err2.2}
For any $\tau>0$
                \begin{equation}
        \label{err2.16}
 \lim_{\omega\to 0} \sup_{t\le \omega^{-1}\tau}\big|x_\omega(t)-X_+(h(\o t))\big|=0
     \end{equation}

        \end{thm}

\begin{figure}[!htbp]
\centering
\includegraphics[width=100mm]{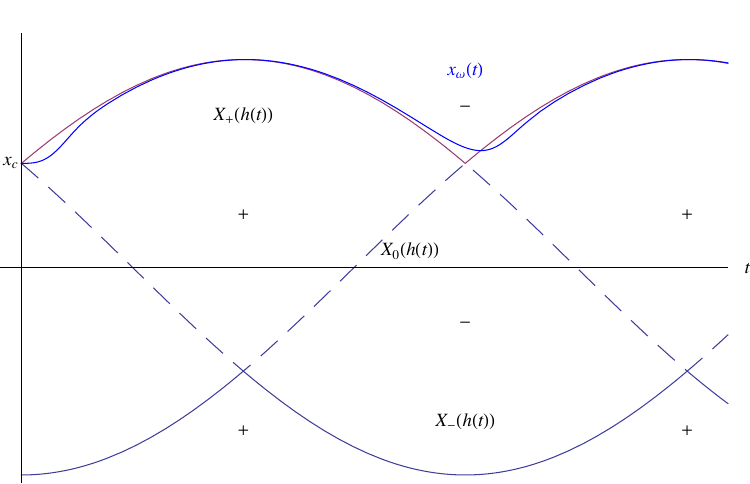}% "%" necessario
\caption{The function $x_\omega(t)$ (blue line) for small
values of $\omega$ tracks the positive branch $X_+(h(t))$ (black line).} \label{m3}
\end{figure}

 Theorem \ref{thm:err2.2} proves that, for oscillations of critical amplitude $h_c$, in the adiabatic limit $\omega\to
 0$ there is not hysteresis
(see Figure \ref{m3}). The relevant time scale is $t=\omega^{-1}
\tau$ and the limit evolution is
                \begin{equation}
        \label{err2.17}
 \lim_{\omega\to 0}x_\omega(\omega^{-1}\tau)=X_+(h(\tau))
     \end{equation}

{\bf The main theorem.}  Theorem \ref{thm:0} asserts that the
dynamics in the macroscopic limit $N\to \infty$  on finite time intervals is described by
the deterministic mean field evolution equation \eqref{eq:eqF}. When $\o$ is small with $N$, the period of the magnetic field oscillations
 is large with $N$. Therefore the behavior exhibited by
\eqref{err2.17} in the adiabatic limit may not correspond to what
the Glauber process does for large but finite $N$. As it will turn
out, it all depends on the way $\omega\to 0$ as $N\to \infty$. As
stated in the introduction, the critical case is $\omega=N^{-2/3}$
to which we restrict hereafter (the origin of the factor $2/3$
will become clear from the proofs but it will also be explained in
Section \ref{section:scalings} in a heuristic way).

There are criticalities   for values of the magnetic field in a
neighborhood of $\pm h_c$. Since $h$ is a periodic function of
time and the process is invariant under change of sign we shall
restrict ourselves to study the behavior in a semi-period. We
consider $t\in N^{2/3}[-\frac \pi 2,\frac \pi 2]$ and suppose
$h=h_N(t)$, with
                \begin{equation}
        \label{err2.18}
 h_N(t):=h(N^{-2/3}t)=-h_c \cos\big(N^{-2/3}t)
     \end{equation}
so that the critical time is  set at $t=0$. We shall denote by
$\p_N$ the law of the process $m_N(t),t\in N^{2/3}[-\frac \pi
2,\frac \pi 2]$ of generator $\L_{h_N(t)}$,  with
$m_N(-N^{2/3}\frac \pi 2)=m_N^0$. We choose such initial value
 in a neighborhood of size $N^{-1/2+\ga}$, $\ga>0$, of the
positive branch, i.e. $|m_N^0-X_+(0)|\le N^{-1/2+\ga}$ (since
$h_N(-N^{2/3}\frac \pi 2)=0$). The main result is given by the
following Theorem. It provides the probability, for large $N$, to
find the magnetization in a neighborhood of one of the two
equilibrium branches $X_\pm(h_N(t))$, respectively, before and
after the critical time $t=0$.

 \vskip.6cm

         \begin{thm}
         \label{thm:TEO} [Main theorem]
       Consider the events
\begin{equation}\label{H+-}
\mathcal H_{\ga}^\pm(I):=\bigg\{ \sup_{t \in I}| m_N(t)-X_{\pm}(
h_N(t))| \le N^{-1/2+\ga}\bigg\},\quad I \subseteq \mathbb{R},
\quad \ga>0
\end{equation}
There is $p_-\in (0,1)$ so that for any $\ga,\eta>0$ and
$\ga'>\ga$, if $|m_N^0-X_+(0)|\le N^{-1/2+\ga}$ then
           \begin{equation}
           \label{eq:TEO1}
\lim_{N \to \i}\p_N\left\{
 \mathcal H^+_{\ga'}\(N^{2/3}\[-\frac \pi 2,-\eta\]\) \right\}=1
          \end{equation}
          \begin{equation}
          \label{eq:TEO2}
\lim_{N \to \i}\p_N\left\{
 \mathcal H^\pm_\ga\(N^{2/3}\[\eta,\frac \pi 2\]\)\right\}=p_{\pm},
          \end{equation}
where $p_+=1-p_-$.
          \end{thm}

The critical interval is $N^{2/3}(-\eta,\eta )$, $\eta>0$
arbitrarily small.  \eqref{eq:TEO1} shows that, in the limit as $N
\to \i$, the magnetization remains, almost surely, in a
neighborhood of size $N^{-1/2+\ga'}$, $\ga'>\ga$, of the positive
branch before the criticality (i.e. for $t<-\eta N^{2/3}$).
\eqref{eq:TEO2} provides the behavior after the criticality (for
$t>\eta N^{2/3}$), it states that  there exists a non-trivial
probability to find the magnetization either in the positive or in
the negative equilibrium branches.

The result can be iterated, as the same arguments can be repeated
every time  the process runs into a criticality. The macroscopic
dynamics is no more deterministic since, at every step there is a
positive probability for the magnetization to jump or not, and the
hysteresis loops observed become, in this sense, random.

\section{Outline of proof}
 \label{section:scalings}

The proof of \eqref{eq:TEO1} is simple.  Indeed, if we fix
$h>-h_c+\eps$, for some $\eps>0$, and the magnetization is
initially in a neighborhood of $X_+(h)$, then $m_N(t)$ has a drift
towards $X_+(h)$. Therefore, with large probability, it stays in a
neighborhood of size $N^{-1/2+\ga'}$ (as $N^{-1/2}$ is the
strength of the noise) of the positive branch. Only after a longer
(exponential) time, tunneling to  the negative branch will be
observed.  In our case $h$ is not fixed but it is so slowly
varying that the above argument remains valid as long as
$h(t)>-h_c+\eps$, for some $\eps>0$ (see Section \ref{app:D}).
When $h$ approaches $-h_c$ the above picture is wrong because at
$-h_c$ the value $x_c$ is stationary but not stable. Lack of
stability and slow changes of the frequency make the noise
competitive with the drift (for the special choice
$\omega=N^{-2/3}$) as we are going to see.

\vskip.6cm

{\bf Scalings.} In order to understand the scalings let us go back
to the stochastic ODE \eqref{err1.1}. Let the magnetic field
oscillate as $h(\o t)=-h_c \cos(\omega t)$, by expanding to
leading orders $F(x)+h$ ($F(x)$ given in \eqref{eq:F}) around $x_c,-h_c$ (i.e.\ for $x-x_c$ and
$\omega t$ both small) we get approximately
         \begin{equation}
         \label{err3.2}
dx= \{h_c\,\frac {(\omega t)^2}2 +\frac{F''(x_c)}2 \,(x-x_c)^2\}
dt + N^{-1/2} dw(t)
         \end{equation}
We scale $y= \omega^{a}(x-x_c)$ and $\tau=\omega^{b}t$,  thus
                  \begin{equation}
         \label{err3.3}
\omega^{-a} dy= \{\omega^{-2b}\frac{h_c \omega^2 \tau^2}2 +
\omega^{-2a}\frac{F''(x_c)y^2}2 \} \omega^{-b}d\tau +
N^{-1/2}\omega^{-b/2} dw(\tau)
         \end{equation}
which becomes independent of $\omega$ and $N$ if
                  \begin{equation}
         \label{err3.4}
         \omega^{a-b/2} N^{-1/2}=1,\quad 2+a-3b=0,\quad a+b=0
         \end{equation}
which yields $\omega=N^{-2/3}$.

The same scalings apply to our case as we shall prove using
extensively martingales techniques.  In order to get rid of
constants in the final equation,  it is convenient to introduce
suitable coefficients in the scaling transformation
\eqref{err3.3}, we define, thus, the process
                  \begin{equation}
         \label{err3.5}
Y_N(t)=\nu N^{1/3} \( m_N\big(\mu N^{1/3} t\big)-x_c\)
         \end{equation}
         with
\begin{equation}\label{mu} \mu=
\(\frac{2}{\b h_c x_c}\)^{1/4}\quad \quad\text{and}\quad
\quad\nu=(\b x_c)^{3/4} \(\frac{2}{h_c}\)^{1/4}
         \end{equation}

We shall study the process $Y_N(t)$ in a time interval which
starts from time $-T$, letting $T\to +\infty$ after $N\to \infty$.
The proof of \eqref{eq:TEO1} can be extended (see Section
\ref{app:D}) till time $-\mu TN^{1/3}$ (which is the microscopic
time corresponding to time $-T$ for $Y_N(\cdot)$) in the following
sense:

\vskip.6cm

           \begin{thm}
           \label{thm:A19}
 There is $c>0$  so that, for any $T$ large enough, $\eps>0$ small
 enough,
         \begin{equation}
         \label{eq:Y}
\limsup_{N \to \i}\p_{N}\Big\{|Y_N(-T)-  T|\le \eps\Big\}\ge
1-e^{-c \eps^2 T}
         \end{equation}

               \end{thm}

 One of the main points in  the proof of \eqref{eq:TEO2} will be to show
(see Sections \ref{section:LD} and \ref{section:Concl}) that the
law of $Y_N(t)$ converges, as $N\to \infty$, to the law of the
stochastic ODE
            \begin{equation}
            \label{eq:xx}
dY(t)=[t^2-Y^2(t)]dt+\xi dw_t, \quad \quad \xi=\frac 2 \b \,
\mu\nu^2
            \end{equation}
which is (modulo multiplicative coefficients) the same as
\eqref{err3.2} with parameters as in \eqref{err3.4}. Due to the
quadratic dependence on $Y$ the solution can blow up in a finite
time, therefore the process is defined with values on $\mathbb
R\cup\{-\infty\}$, with the convention that, if $Y(t)=-\infty$,
then $Y(t')=-\infty$ for all $t'\ge t$.  The drift in
\eqref{eq:xx} vanishes on the two straight lines $Y=\pm t$. It is
negative for $Y<-|t|$ and it points towards  $|t|$ for $Y>-|t|$. A
more careful analysis shows that there is a critical trajectory
$y^*(t)<0$ solution of the deterministic version (i.e.\ with $\xi
= 0$) of \eqref{eq:xx} such that any deterministic solution which
starts above the critical curve  is exponentially asymptotic to
$(t,t)$ as $t\to \infty$.

%\vskip.3cm
%
%\begin{figure}[htbp]
%\centering
%\includegraphics[width=100mm]{YGra2}% "%" necessario
%\caption{Phase diagram of the solutions $y(t)$ of the
%deterministic version (i.e. the $\xi=0$ case) of \eqref{eq:xx}
%(see \eqref{eq:y}). Notice the critical solutions $y^*(t)$
%($|y^*(t)+t|\to 0$ as $t \to +\i$) and $y^+(t)$ ($|y^+(t)+t|\to 0$
%as $t \to -\i$).}\label{YGra2}
%\end{figure}

We denote by $\p_{-T,y}$ the law on $\mathbb R\cup\{-\infty\}$ of
the solution $Y(t),\, t>-T$ of \eqref{eq:xx} starting from
$Y(-T)=y$, $T>0$.  In Section \ref{section:Asy} we prove the
following Theorem.

\vskip.6cm

           \begin{thm}
           \label{thm:3.2}
Let  $\p$ be the probability law with support on solutions $Y(t)$
of \eqref{eq:xx} such that
       \begin{equation}\label{eq:3.14}
\lim_{t\to -\infty}|Y(t)+t| =0 \quad \quad \p-\text{a.s.}
 \end{equation}
 then there exist $p_\pm\in (0,1)$, $p_+=1-p_-$, such that
 \begin{equation}
\p \Big\{ \; \text{there is $t:\;Y(t) =-\infty$}\; \Big\}=p_-
\quad \text{and} \quad \p \Big\{\lim_{t\to
\infty}|Y(t)-t|=0\Big\}=p_+
            \end{equation}
For any $\eps>0$ small enough, for any bounded continuous function
$g(y)$ with compact support and any $t\in \mathbb R$,
         \begin{equation}
         \label{eq:Y?}
\lim_{T \to \i} \mathbf 1_{|y-T|\le \eps}
\,\E_{\p_{-T,y}}\[g(Y(t))\]= \E_{\p}\[g(Y(t))\]
         \end{equation}
Moreover  there exists $c>0$ such that, for any $T$ large enough,
$\eps$ small enough,
\begin{equation}
         \label{eq:Y??}
\p \Big\{|Y(-T)-T|\le \eps\Big\}>1-e^{-c\eps^2T}
         \end{equation}

               \end{thm}

Thus with $\p$ probability one either $Y(t)$ blows up in a finite
time or it is asymptotic to $t$ as $t\to \infty$, both events
having non zero probability. The next goal is to extend the above
result to the finite $N$ process $Y_N(t)$. For $T>0$ we define the
rectangle:
            \begin{equation}
            \label{eq:xx0}
\mathcal R_{T}=\Big\{(t,y)\in \mathbb R^2: t\in [-T,T], |y|\le
2T\Big\}
            \end{equation}

and,  for $\eps \in(0,1)$,
\begin{eqnarray}
            \label{eq:xx00}
\partial \mathcal R^+_{T}:=\{T\} \times [T-\eps, T+\eps] \quad
\partial \mathcal R^-_{T}:=[-T,T]\times \{-2T\},
            \end{eqnarray}
$\partial \mathcal R_T^\pm \subseteq \partial \mathcal R_T$. For
the processes $Y(t)$ such that $(-T,Y(-T))\in \mathcal R_T$, we
denote by $\tau_{T}$ the first exit time from $\mathcal R_{T}$
\begin{equation}\label{tau}
\tau_T:=\inf \left\{t\ge -T: \: Y(t)\notin \mathcal R_T\right\}
\end{equation}
 and
define the sets
            \begin{equation}
            \label{eq:xx?}
\mathcal E^\pm_{T}=\Big\{Y\; : \:   (\tau_{T},Y(\tau_{T}))\in
\partial \mathcal R^\pm_{T}\Big\}
            \end{equation}

 We shall prove in Section \ref{section:Asy}

\vskip.6cm

            \begin{prop}
            \label{prop:3.1}
Let $Y(t),\, t\ge -T$ be a solution of \eqref{eq:xx} starting at
$-T$ from $y:|y-T|\le \eps$, $\eps>0$ small enough, then
             \begin{equation}
            \label{eq:xx??}
 \lim_{T \to \i}\p_{-T,y}\Big\{Y \in \mathcal E^+_{T}\cup \mathcal
 E^-_{T}\Big\}=1
            \end{equation}
moreover
             \begin{eqnarray}
            \label{eq:xx??+}
&& \hskip.7cm\lim_{T \to \i}\p_{-T,y}\Big\{\lim_{t \to
\i}|Y(t)-t|=0\; \Big|\: \,Y \in \mathcal
E_T^+\Big\}=1\\
&& \text{and} \quad  \lim_{T \to \i}\p_{-T,y}\Big\{ \; \text{there
is $t:Y(t) =-\infty$}\; \Big| \: \,Y \in \mathcal E_T^-\Big\}=1.
            \end{eqnarray}
             \end{prop}

The following Corollary is a direct consequence of Theorem
\ref{thm:3.2} and Proposition \ref{prop:3.1}:

\vskip.6cm

            \begin{coro}
            \label{prop:3.4}
For $Y(t)$ as in the previous Proposition we have
                        \begin{equation}
            \label{eq:3.15}
\lim_{T \to \i} \big| \p_{-T,y}\{ Y \in \mathcal
E^{\pm}_{T}\}-p_{\pm}\big| =0.
            \end{equation}

             \end{coro}

Let $\p_{N,-T,y}$ be the law of $Y_N(t)$ given $Y_N(-T)=y$. Using
martingale convergence theorems, in Section \ref{section:Concl} we
prove the following result.

\vskip.6cm

            \begin{prop}
            \label{prop:3.5}
For $Y(t)$ solution of \eqref{eq:xx} starting  at $-T$ from
$y:|y-T|\le \eps$, $\eps$ small enough, we have
                        \begin{equation}
            \label{eq:3.16}
\lim_{N\to \infty} \p _{N,-T,y}\{Y_N \in \mathcal E^{\pm}_{T}\}=
\p _{-T,y}\{Y \in \mathcal E^{\pm}_{T}\}
            \end{equation}
            and
             \begin{equation}
            \label{eq:3.16'}
\lim_{N\to \infty} \p _{N,-T,y}\{Y_N \in \mathcal E^{+}_{T}\cup
\mathcal E^-_T\}= \p _{-T,y}\{Y \in \mathcal E^{+}_{T}\cup
\mathcal E^-_T\}.
            \end{equation}

             \end{prop}

Proposition \ref{prop:3.5} allows us to extend the results
obtained for $Y(t)$ to the finite $N$ process $Y_N(t)$. Finally in
Section \ref{section:Concl} we prove the following Proposition
that is the last ingredient to conclude the proof of Theorem
\ref{thm:TEO}.

\vskip.6cm

            \begin{prop}
            \label{prop:3.6}
For any $\eta, \ga>0$,
             \begin{equation}
            \label{m}
 \lim_{T \to \i}\lim_{N \to \i}\p_{N}\left\{\mathcal H^\pm_\ga\(N^{2/3}\[\eta,\frac \pi 2\]\)  \,\bigg| \:Y_N \in \mathcal
 E^\pm_{T}\right\}=1.
            \end{equation}

             \end{prop}

\section{Limit dynamics in the critical region}
\label{section:LD}

The study of the limit behavior as $N \ra \i$ of the spin-flip
evolution defined in Section 2 is based on some martingale
theorems. In our dynamics we have two natural martingales:
\begin{equation}\label{eq:ML}
M_{N,T}(t)=m_N(t)-m_N(-\mu T N^{1/3})-\int_{-\mu T N^{1/3} }^t
\F_N(m_N(s),\h(s))\; ds
\end{equation}
where $\F_N (x,h):=\L_h x$, $T >1$, and
\begin{equation}\label{eq:mart2}
M_{N,T}^2(t)-\int_{-\mu T N^{1/3}}^t \mathcal G_N(m_N(s),\h(s))\;
ds
\end{equation}
with $\mathcal G_N(x,h):=\L_h x^2-2x \L_h x$.

\vskip.5cm
In the following Lemma we prove that, for large $N$, the function $\mathcal F_N(x,h)$ is well approximated
by the infinite volume drift $F(x,h)=-x+\tanh\{\b(x+h)\}$ (see the infinite volume equation \eqref{eq:F}).

\vskip.6cm
\begin{lemma}\label{lemma:F}
There exists $c>0$ such that, for any $x \in[-1,1]$, $|h|\le h_c$,
$N$ large enough,
\begin{equation}\label{eq:1}
\big|\F_N(x,h)-F(x,h)\big| \le \frac c N
\end{equation}
and, for $\Lambda(x,h)=1-x\tanh\{ \b(x+h)\}$,
\begin{equation}\label{eq:4}
\big|N \mathcal G_N(x,h)-2\Lambda(x,h)\big| \le \frac {c} N
\end{equation}

\end{lemma}

 {\bf Proof.} We have
\begin{eqnarray*}
\F_N(x,h)=\frac 2 N \( c^+(x,h)-c^-(x,h) \)=  \( \hat
c^+(x,h)-\hat c^-(x,h) \)-x\( \hat c^+(x,h)+\hat c^-(x,h) \)
\end{eqnarray*}
then there exists $c>0$ such that
\begin{equation}\label{eq:1(1)}
\big|\hat c^+(x,h)+\hat c^-(x,h)\big|\le \frac c N \quad\quad
\text{and} \quad \quad \big|(\hat c^+(x,h)-\hat c^-(x,h))- \tanh\{
\b(x+h)\}\big|\le \frac c N
\end{equation}
for any $N$ large enough, that yields $\eqref{eq:1}$. Now
\begin{equation*}
\mathcal G_N(x,h)=\frac 4 {N^2}[c^+(x,h)+c^-(x,h)]
\end{equation*}
thus
\begin{equation*}
2-N \mathcal G_N(x,h)=2\[ 1-(\hat c^+(x,h)+\hat c^-(x,h)) \] +2x
\[ \hat c^+(x,h)-\hat c^-(x,h) \]
\end{equation*}
then $\eqref{eq:4}$ follows from $\eqref{eq:1(1)}$. \qed

\vskip.3cm

Let $Y(t), t \ge -T$ be the solution of \eqref{eq:xx} starting
from $Y(-T)=y$,  and $\tau_{T}$ be the first exit time from the rectangle $\mathcal R_T$ (see
$\eqref{tau}$ and \eqref{eq:xx0}). We denote by $\p^*_{-T,y}$ the
law of the stopped process $Y(t \wedge \tau_{T})$ on $\mathcal
D[-T,T]$. We call $\tau_{N,T}$ the corresponding stopping time for
the finite $N$-process $Y_N(t)$ (see \eqref{err3.5}) and denote by
$\p^*_{N,-T,y}$ the law of the corresponding stopped process. We
are going to prove (see Proposition \ref{prop:z,x}) the
convergence of $\p^*_{N,-T,y}$ to $\p^*_{-T,y}$ for  suitable
 $T,y$. Let $\mathcal D[-T,T]$ be the space of functions on $[-T,T]$ that are right-continuous and have left-hand
  limits. The convergence results in this Section are meant
in the sense of the Skorohod metric on $\mathcal D[-T,T]$.
For more details on the space $\mathcal D[-T,T]$ and the weak convergence on
$\mathcal D[-T,T]$ see Chapter 3 of \cite{Bi}.

\vskip.3cm

For  the martingale $\hat M_{N,T}(t):=\nu N^{1/3}M_{N,T}(\mu
N^{1/3}(t\wedge \tau_{N,T}))$,
\begin{eqnarray*}
\hat M_{N,T}(t)=Y_N(t \wedge \tau_{N,T})-Y_N(-T)-\nu \mu
N^{2/3}\int_{-T}^{t \wedge \tau_{N,T}} \mathcal F_N(m(\mu
sN^{1/3}),h_N(\mu sN^{1/3}))\, ds
\end{eqnarray*}
we have the following result

\vskip.6cm

\begin{prop}\label{prop:1}

Let $w(t)$ be the standard Brownian motion and $\xi:= \frac 2 \b
\mu \nu^2$, then
\begin{equation}\label{eq:Mp}
\begin{CD}
\hat M_{N,T}(t)   @>\mathcal D>>\xi \, w(T+ t\wedge  \tau_{T}) \quad  \text{as}\quad N \ra \i.\\
\end{CD}
\end{equation}

\end{prop}

{\bf Proof.} By $\eqref{eq:mart2}$, the quadratic variation of
$\hat M_{N,T}(t)$ is given by
\begin{eqnarray*}
\hat V_{N,T}(t):=\nu^2 \mu N \int_{-T}^{t \wedge
\tau_{N,T}}\mathcal G_N(m_N(\mu sN^{1/3} ),h_N( \mu sN^{1/3} ))\,
ds
\end{eqnarray*}
thus, for  $\Lambda(m,h)$ as  in Lemma $\ref{lemma:F}$,  by
$\eqref{eq:4}$, there exists $c>0$ such that
\begin{eqnarray*}
\sup_{t \ge -T}\bigg| \hat V_{N,T}(t)-2\nu^2 \mu \int_{-T}^{t
\wedge \tau_{N,T}} \Lambda(m_N( \mu sN^{1/3} ),h_N(\mu s
N^{1/3}))\;ds\bigg| \le cN^{-1}
\end{eqnarray*}
for any $N$ large enough. In a neighborhood of $(x_c,-h_c)$,
\begin{equation*}
\Lambda(x,h)= \frac 1 \b + \OO(h+h_c)+\OO(x-x_c)
\end{equation*}
moreover, for $t<N^{2/3}$ there exists $c>0$ such that
$|h_N(t)+h_c|\le c (tN^{-2/3})^2$ for any $N$ large enough. We
have, thus
\begin{equation*}
\sup_{-T \le s \le \tau_{N,T}}\bigg |\Lambda(m_N(\mu s
N^{1/3}),h_N(\mu s N^{1/3}))- \frac 1 \b \bigg|\le c N^{-1/3}
\end{equation*}
for a suitable $c>0$ then
\begin{equation}\label{eq:D2}
\sup_{t \ge -S}\big|\hat V_{N,T}(t)-\xi(T+t \wedge \ta)\big| \le
cN^{-1/3}
\end{equation}
We have $\ta \stackrel{\mathcal P}{\ra}\tau_{T}$,  for $N\to \i$,
hence, by $\eqref{eq:D2}$,
\begin{equation*}
\hat V_{N,T}(t) \stackrel{\mathcal P}{\ra} \xi (T+t \wedge
\tau_{T}) \quad \text{as} \quad N \to \i
\end{equation*}
thus $\eqref{eq:Mp}$ follows   since  $\hat M_{N,T}(-T)=0$ and
$\hat M_{N,T}(t)$ has at most discontinuities of order $N^{-2/3}$
\\(see \cite{Bi} and \cite{Po}). \qed

\vskip.6cm

\begin{prop}\label{prop:z,x}
For any $T,y>0$ such that  $y<2T$,  $\p^*_{N,-T,y}$ converges to
$\p^*_{-T,y}$ as $N \to \i$.
\end{prop}

{\bf Proof.} As usual with martingale problems, we first need to
prove tightness and then to identify the limiting points by
proving that they satisfy a martingale equation which has unique
solution. By Proposition \ref{prop:1} follows the tightness of
$\hat M_{N,T}(t)$. It remains to prove the tightness of
\begin{equation*}
\Gamma_{N,T}(t)=\nu \mu N^{2/3}\int_{-T}^{t \wedge \tau_{N,T}}
\mathcal F_N(m_N(\mu sN^{1/3}),h_N(\mu sN^{1/3}))\, ds
\end{equation*}
We use the Chensov moment condition, indeed there exists  $c$ such
that, for all $t>s\ge -T$,
\begin{equation}\label{eq:E}
\E_{\p_{N,-T,y}} \[  \left|
\Gamma_{N,T}(t)-\Gamma_{N,T}(s)\right|^2
\] \le c |t-s|^2
\end{equation}
where $\eqref{eq:E}$ holds after using the Cauchy-Schwartz
inequality, being the integrated function in $L^2$.
 It follows that the  stopped process $Y_N(t \wedge\ta)$ is tight and,
  consequently,  its law $\p_{N,-T,y}$  converges by subsequences.
Moreover, any limiting point has support on $\mathcal
C([-T,T],\mathbb{R})$, this follows from the fact that the jumps
of $Y_N$ are $\pm N^{-\frac 2 3}$.

By $\eqref{eq:1}$, we can approximate the term $\mathcal F_N(x,h)$
in $\eqref{eq:ML}$ with $F(x,h)$ unless errors of order $N^{-1}$.
We perform the Taylor expansion of $F(x,h)$ in a neighborhood of
$(x_c,-h_c)$. Being $F(x_c, -h_c)=\partial F/ \partial x
(x_c,-h_c)=0$, the leading terms are the first order in $(h+h_c)$
and the second order in $(x-x_c)$, we have
\begin{eqnarray*}
F(x,h)=(h+h_c)-\b x_c (x-x_c)^2+\mathcal
O((h+h_c)(x-x_c))+\mathcal O((h+h_c)^2)+ \mathcal O((x-x_c)^3)
\end{eqnarray*}
On the other hand, for $t N^{-2/3}$ vanishingly small as $N \to
\i$, $h_N(t)=-h_c +{h_c}t^2 N^{-4/3}/2 +\mathcal O((t
N^{-2/3})^4)$, thus there exists $c$ such that
\begin{equation}\label{eq:F1}
\sup_{t \in \mu N^{1/3}[-T,\tau_{N,T}]}\bigg|\mathcal
F_N(m_N(t),h_N(t))-\bigg\{\frac {h_c}2 \, t^2 N^{-4/3}-\b x_c
(m_N(t)-x_c)^2\bigg\}\bigg|\le c N^{-1}
\end{equation}
for $N$ large enough, then, by $\eqref{eq:F1}$,
\begin{equation}\label{eq:F2}
\sup_{t \ge -T}\bigg|\hat M_{N,T}(t)-Y_N(t \wedge
\tau_{N,T})+Y_N(-T) +\int_{-T}^{t \wedge \tau_{N,T}} \bigg\{ \frac
{h_c}2 \, \nu \mu^3 s^2- \b x_c \, \mu \nu^{-1} Y_N^2(s)\bigg\}\,
ds\bigg| \le c N^{-1/3}
\end{equation}
For our choice of $\mu$ and $\nu$ (see \eqref{err3.5}), the
integrand in $\eqref{eq:F2}$ becomes $s^2-Y_N^2(s)$. From
$\eqref{eq:F2}$ and Proposition $\ref{prop:1}$ we  deduce that any
limiting point satisfies a martingale relation that uniquely
defines a process which is the law of the solution of
$\eqref{eq:xx}$. \qed

\section{Behavior of the limit process}
\label{section:Asy}

In this Section we are going to investigate the behavior of a
generic solution $Y(t)$  of the SDE
\begin{equation}\label{x}
dY(t)=[t^2-Y^2(t)]dt+\xi dw_t, \quad \quad \xi>0,
\end{equation}
For any fixed $t_0 \in \mathbb{R} \cup \{-\i\}$, $y_0\in
\mathbb{R}$, we denote by $\p_{t_0,y_0}$ the probability law of
the process $Y(t), \;t \ge t_0$ solution of \eqref{x} starting
from $y_0$ at time $t_0$. Moreover we denote by $\p$ the law of
$Y(t), \;t \in \mathbb{R}$ solution of \eqref{x} conditioned to
$|Y(t)+t|\to 0$ as $t \to -\i$.

\subsection*{Deterministic analysis}
   %\label{app:A}

One of the preliminary steps for the study of \eqref{eq:xx} is the
analysis of the related deterministic equation
\begin{equation}\label{eq:y}
 y'(t)=t^2-y^2(t)
\end{equation}
Proposition $\ref{prop:y}$ is proved in Section 2.3 of \cite{Ca},
it concerns the asymptotic behavior for $t \to \i$ of a generic
solution $y(t)$ of \ref{eq:y}.

\vskip.3cm

\begin{prop}\label{prop:y}
There exists a  decreasing solution $y^*(t)$ of \eqref{eq:y} such
that $-t>y^*(t)>-\sqrt{t^2+1}$, for any $t\ge 0$.   Let $y(t)$ be
the solution of \eqref{eq:y} starting at time $t_0\ge 0$ from $y_0
\in \mathbb{R}$,
\begin{itemize}
  \item if $y_0>y^*(t_0)$, then, for any $\d \in (0,1)$ there exists  $t_\d \ge t_0$ such that $|y(t)-t|\le \frac 1
  {2(1-\d)t}$ for any $t\ge t_\d$;
  \item if $y_0<y^*(t_0)$, then $y(t)$ is decreasing for $t\ge 0$ and it explodes to $-\i$ in a
  finite time.
\end{itemize}

\end{prop}

\vskip.3cm

\subsection*{Asymptotic behavior of $Y(t)$ for $t \to \i$}

In this first part of the Section we prove the following Theorem.

\vskip.3cm
\begin{thm}\label{thm:stoc}
Consider the sets
\begin{eqnarray}\label{eq:E+-}
E^+:=\{Y\::\:\lim_{t \to \i} |Y(t)- t|=0\} \quad \quad
\text{and}\quad \quad E^-:=\{Y\::\; \text{there is $ \: t:Y(t)
=-\infty$}\;\}
\end{eqnarray}
 then $\p_{t_0,y_0}\{Y \in E^+ \cup E^-\}=1$ for any $t_0 \in \mathbb{R}\cup\{-\i\}$, $y_0 \in
 \mathbb{R}$.
\end{thm}

 The proof of Theorem $\ref{thm:stoc}$ consists of three parts.
We define the stopping time
\begin{equation*}
\Pi:= \inf\{t \: : \: Y(t)=-\i\}
\end{equation*}
then $\Pi \in \mathbb{R}\cup \{+\i\}$. We fix $T>0$ large enough,
suppose $\Pi>T$ and study the behavior of $Y(t)$ for $t\ge T$. In
Proposition \ref{prop:I} we prove that if $Y(t)$ is in a
 neighborhood of $y^*(t)$ at time $T$ then
$Y(t)$ escapes from it $\p$- a.s. In Propositions \ref{prop:II}
and \ref{prop:III} we prove that the probability for the events $Y
\in E^\mp$ to occur is close to the probability that $Y(t)$ leaves
such a critical neighborhood, respectively, from below or from
above. Unless further indications, in this Section we mean, by
$c$, a positive constant not depending on $T$.

\vskip.3cm

We will denote by $y^*(t)$ the solution of the ODE \eqref{eq:y}
defined in Proposition \ref{prop:y}, and define the processe
$z^*(t):=Y(t)-y^*(t)$. $z^*(t)$ verifies the equation
\begin{eqnarray}
dz^*(t)=-z^*(t)(z^*(t)t+2y^*(t))dt+\xi dw(t).\label{eq:z'}
\end{eqnarray}
For any fixed $\d>0$ small enough, we define the stopping time
$\tau^*_{T,\d}:=\inf \left\{t\ge T \; : \: |z^*(t)|\ge \d
\right\}$.

\vskip.6cm

\begin{prop}\label{prop:I}
For any $T>0$,  $\d>0$ small enough,
\begin{equation}\label{step1}
\mathbf 1_{\Pi>T}\;\p_{T,Y(T)} \Big\{\tau^*_{T,\d}<\i\Big\}=1
\end{equation}
\end{prop}

{\bf Proof.} Let us assume $\Pi>T$.  We need to prove the
assertion for the paths such that $|z^*(T)|<\d$. Suitably applying
the Ito's formula to \eqref{eq:z'}, we get
\begin{equation}\label{eq:z2}
dz^{*2}(t)=[-2z^{*2}(t)(z^{*2}(t)+2y^*(t))+\xi^2]dt+2\xi z^*(t)
dw(t),
\end{equation}
thus, for $T \le t \le \tau^*_{T,\d}$
\begin{eqnarray}\label{eq:Asy11(1)}
&&z^{*2}(t)
 \ge z^{*2}(T)+\xi^2 t +2 \xi \int_{T}^{t} z^*_s dw_s,
\end{eqnarray}
the  inequality descending since, for $\d$ small enough,
$-2z^{*2}_{t \wedge \tau^*_{T,\d}}(z^*_{t \wedge \tau^*_{T,\d}}+2
y^*_{t \wedge \tau^*_{T,\d}})\ge 0$.

The process $2 \xi \int_T^{t\wedge \tau^*_{T,\d}} z^*_s dw_s$ is a
continuous martingale, thus  its expected value is constantly zero
and
\begin{equation*}
\E\[\bigg(2 \xi \int_T^{t\wedge \tau^*_{T,\d}} z^*_s
dw_s\bigg)^2\]=4 \xi^2 \int_T^t \E\[z^{*2}_s \mathbf 1_{s \le
\tau^*_{T,\d}}\]\, ds \le 4 \xi^2 \d^{2} (t-T)
\end{equation*}

hence, by the Doob's inequality, for any $n \in \mathbb{N}$,
\begin{equation*}
\p_{T,Y(T)}\bigg\{2 \xi \, \bigg|\int_T^{(T+n^4)\wedge
\tau^*_{T,\d}} z^*_s dw_s \bigg|\ge n^3\bigg\}\le \frac{4\xi^2
\d^2}{n^2}
\end{equation*}

thus, from the Borel-Cantelli Lemma and $\eqref{eq:Asy11(1)}$,
$\p_{T,Y(T)}$-a.s., there exists $\tilde n$ such that, for $n \ge
\tilde n$,
\begin{equation*}
\d^2 \ge z^{*2}((T+n^4) \wedge \tau^*_{T,\d})
>-n^3+\xi^2((T+n^4)\wedge \tau^*_{T,\d})
\end{equation*}
then $\tau^*_{T,\d} \le (T+n^4) \vee(\d^2 + n^3)/\xi^2$, thus, for
any $T>0$
\begin{eqnarray*}
\p_{T,Y(T)}\Big\{\tau^*_{T,\d} <\i\Big\} \ge
\p_{T,Y(T)}\Big\{\liminf_{n \to +\i} \;\{\tau^*_{T,\d} \le
(T+n^4)\}\Big\}=1
\end{eqnarray*}
and \eqref{step1} is proved. \qed

\vskip.6cm

We omit the proof of the following Lemma.

\vskip.3cm

\begin{lemma}\label{lemma:exp}
Let $t>s$, for any $\gamma>0$, we have
\begin{equation}\label{eq:exp+}
\frac {e^{\gamma t^2}}{2\gamma t}\,[1-e^{-\gamma (t^2-s^2)/2}]\le
\int_{s}^{t}e^{\gamma u^2}\, du \le \frac {e^{\gamma t^2}}{2\gamma
t}\(\frac{2\gamma s^2}{2\gamma s^2-1}\)
\end{equation}
for any $s>\frac 1 {\sqrt {2 \gamma}}$, and
\begin{equation}\label{eq:exp-}
\frac {e^{-\gamma s^2}}{2\gamma s}[1-t^{-1}e^{-\gamma
(t^2-s^2)/2}]\(\frac{2\gamma s^2}{2\gamma s^2+1}\)\le
\int_{s}^{t}e^{-\gamma u^2}\, du \le \frac{e^{-\gamma
s^2}}{2\gamma s}
\end{equation}
for any $s>0$.
\end{lemma}

\vskip.6cm

\begin{prop}\label{prop:II}
There exists $c>0$ such that, for any   $T$ large enough, $\d>0$,
\begin{equation}\label{step2}
\mathbf 1_{\Pi>T}\;\mathbf 1_{\tau^*_{T,\d}<\i, \;
z^*(\tau^*_{T,\d})<-\d}
\:\p_{\tau^*_{T,\d},Y(\tau^*_{T,\d})}\Big\{ Y \notin E^- \Big\}\le
e^{-c T}
\end{equation}
\end{prop}

{\bf Proof.} Suppose $\tau^*_{T,\d}<\i$ and $\Pi>T$, thus,
consequently, $\Pi>\tau^*_{T,\d}$. We denote by $\hat y(t)$  the
 solution of the ODE $\eqref{eq:y}$ starting at time
$\tau^*_{T,\d}$ from $y^*(\tau^*_{T,\d})-\d/2$. From Proposition
\ref{prop:y} we know that $\hat y(t)$ explodes to $-\i$ in a
finite time. Consider $\hat z(t):=Y(t)-\hat y(t)$, thus $\hat
z(t)$ verifies the SDE
\begin{equation*}
d\hat z_t=-\hat z_t(\hat z_t+2\hat y(t)\hat z_t)dt+\xi dw_t
\end{equation*}
 We can assume $\hat z(\tau^*_{T,\d})\le -\d/2$ since $1_{z^*(\tau^*_{T,\d})<-\d}\le \mathbf 1_{{\hat z}(\tau^*_{T,\d})\le
 -\d/2}$. We have
\begin{eqnarray}\label{z}
\hat z(t)=\hat z(\tau^*_{T,\d})\;e^{-2\int_{\tau^*_{T,\d}}^t \hat
y(s)\, ds}-\int_{\tau^*_{T,\d}}^t
\hat z^{2}(u)\;e^{-2\int_{u}^t \hat y(s)\, ds}\; du + \xi \chi_{\tau^*_{T,\d}}(t) \\
\text{with} \quad \quad \hat \chi_{\tau^*_{T,\d}}(t):=
\int_{\tau^*_{T,\d}}^t e^{2 \int_{\tau^*_{T,\d}}^u \hat y(s)\, ds}
\; dw_u\hskip1cm
\end{eqnarray}
then
\begin{equation*}
 \mathbf 1_{\hat z(\tau^*_{T,\d})\le -\d/2}
\: \p_{\tau^*_{T,\d},Y(\tau^*_{T,\d})}\left\{\hat z(t)\le e^{-2
\int_{\tau^*_{T,\d}}^t \hat y(s) \, ds}\( \hat
\chi_{\tau_{T,\d}^*}(t)- \d/ 2 \),\, \forall t \ge
\tau^*_{T,\d}\right\}=1.
\end{equation*}

The probability law of $\hat \chi_{\tau^*_{T,\d}}(t)\, |\;
\tau^*_{T,\d}$, $t \ge \tau^*_{T,\d}$ is a centered gaussian.
Since $\hat y(t)\le -t$, $t \ge \tau^*_{T,\d}$, we have
\begin{equation}\label{Q}
\E\[\hat \chi^2_{\tau^*_{T,\d}}(t)\; \Big| \; \tau^*_{T,\d}
\]=\int_{\tau^*_{T,\d}}^{t} e^{4 \int_{\tau^*_{T,\d}}^t \hat y(s)\,
ds}\; du \le  e^{2 \tau^{*2}_{T,\d}}\int_{\tau^*_{T,\d}}^t
e^{-2u^2} \; du\le \frac 1 {4T}
\end{equation}
 where the last inequality descends from \eqref{eq:exp-}, since $\tau^*_{T,\d}\ge
T$. Hence there exists $c>0$,  such that
\begin{equation}\label{.0}
\p_{\tau^*_{T,\d},Y(\tau^*_{T,\d})}\left\{\sup_{t \ge
\tau^*_{T,\d}}\hat z(t) \ge 0\right\}\le
\p_{\tau^*_{T,\d},Y(\tau^*_{T,\d})}\left\{\sup_{t \ge
\tau_{T,\d}^*}\hat \chi_{\tau_{T,\d}^*}(t) \ge \frac \d 2
\right\}\le e^{-c T}
\end{equation}
for any $T$ large enough, where the second inequality follows from
\eqref{MS} and \eqref{Q}. Then we get \eqref{step2}. \qed

\vskip.6cm

We denote by $y^+(t)$ the solution of  \eqref{eq:y} conditioned to
$\lim_{t \to -\i}|y^+(t)+t|=0$  and define the process
$z^+(t):=Y(t)-y^+(t)$. $z^+(t)$ satisfies the SDE
\begin{eqnarray}
dz^+(t)=-z^+(t)(z^+(t)t+2y^+(t))dt+\xi dw(t),\label{eq:z}
\end{eqnarray}
thus, for any $t_0 \in\mathbb{R} \cup \{-\i\}$,
\begin{eqnarray}\label{zz}
z^+(t)=z^+(t_0)e^{-2\int_{t_0}^t y^+(s)\, ds}-\int_{t_0}^t
z^{+2}(u)e^{-2\int_{u}^t y^+(s)\, ds}\; du + \xi \chi^+_{t_0}(t) \\
\text{with} \quad \quad \chi^+_{t_0}(t):=\int_{t_0}^t
e^{-2\int_u^t y^+(s)\, ds}\; dw_u. \label{zzz}\hskip1cm
\end{eqnarray}

We fix $\eps>0$ small enough and define the stopping time
$\tau^+_{T,\eps}:=\inf \left\{t\ge T \; : \: |z^+(t)|\le \eps
\right\}$.

\vskip.6cm

\begin{lemma}\label{x0}
For $t_0 \in \mathbb{R} \cup \{-\i\}$, $\chi^+_{t_0}(t)$ as
 in \eqref{zzz}, $y_0 \in \mathbb{R}$, there exists $c>0$ such
that, for  $\la$ large enough,
\begin{equation}\label{eq:x0}
\p_{t_0,y_0}\left\{\sup_{t \ge t_0} \:
\chi^+_{t_0}(t)\:(\sqrt{|t|}\vee 1)>\la\right\} \le e^{-c \la^2}
\end{equation}
\end{lemma}

{\bf Proof.} $\chi^+_{t_0}(t)$ is a centered Gaussian process of
variance
\begin{equation*}
\E\[\chi^{+2}_{t_0}(t)\]= \int_{t_0}^t e^{-4\int_u^t y^+(s)\,
ds}\; du.
\end{equation*}
Let us suppose, at first, $t_0<0$. We know that $y^+(t)\ge -t$ for
$t<0$, thus, for $t_0\le t<0$,
\begin{equation*}
\E\[\chi^{+2}_{t_0}(t)\]\le  e^{2 t^2} \int_{|t|}^{|t_0|} e^{-2
u^2}\, du \le \frac 1 {4|t|} \wedge 1
\end{equation*}
where  the second inequality follows from \eqref{eq:exp-}.
 A similar  estimate can be obtained for $t\ge 0$ using
\eqref{eq:exp+}, since $\inf_{t \in \mathbb{R}}y^+(t)>0$ and
$y^+(t)\ge t-1/t$ for $t>0$ large enough. Therefore, for any $t_0
\in \mathbb{R} \cup \{-\i\}$, there exists $c>0$ such that, for
$t\ge t_0$,
\begin{equation*}
\E\[\chi^{+2}_{t_0}(t)\]\le  c\(\frac 1 {|t|} \wedge 1\),
\end{equation*}
thus \eqref{eq:x0} follows from inequality \eqref{MS}. \qed

\vskip.6cm

\begin{prop}\label{prop:III}
There exists $c>0$ such that, for any  $T$ large enough,
$\d,\eps>0$,
\begin{equation}\label{step3}
\mathbf 1_{\Pi>T}\;\mathbf 1_{\tau^*_{T,\d}<\i, \;
z^*(\tau^*_{T,\d})>\d}\:\p_{\tau^*_{T,\d},Y(\tau^*_{T,\d})}\Big\{
\tau^+_{T,\eps}=+\i \Big\}\le e^{-cT}
\end{equation}
\end{prop}

{\bf Proof.} Suppose $\tau^*_{T,\d}<\i$ and $\Pi>T$, then
$\Pi>\tau^*_{T,\d}$. As in the proof of  Proposition \ref{prop:I},
we mainly make use of comparison arguments. We compare, by means
of Lemma \ref{lemma:comp}, the process $z^+(t)$ with suitable
gaussian processes. Then use the inequality \ref{MS} to estimate
the behavior of such gaussian processes. We will avoid the
details, let us see.
 Suppose $z^*(\tau^*_{T,\d})>\d$ and
$|z^+(\tau^*_{T,\d})|>\eps$, we need to distinguish two cases:
$z^+(\tau^*_{T,\d})> \eps$ and $z^+(\tau^*_{T,\d})<-\eps$.

\vskip.1cm

Consider the first case $z^+(\tau^*_{T,\d})>\eps$,  from
\eqref{zz}, we have
\begin{equation*}
1_{z^+(\tau^*_{T,\d})>\eps}\:\p_{\tau^*_{T,\d},Y(\tau^*_{T,\d})}\left\{z^+(t)
\le z^+(\tau^*_{T,\d})\;e^{-2\int_{\tau^*_{T,\d}}^t y^+(s)\, ds}+
\xi \chi^+_{\tau^*_{T,\d}}(t)\right\}=1,
\end{equation*}
thus
\begin{eqnarray}\label{z+}
&&
1_{z^+(\tau^*_{T,\d})>\eps}\:\p_{\tau^*_{T,\d},Y(\tau^*_{T,\d})}\left\{\inf_{t
\ge \tau^*_{T,\d}} z^+(t)>\eps\right\}\nonumber\\&&\le
 1_{z^+(\tau^*_{T,\d})>\eps}\: \p_{\tau^*_{T,\d},Y(\tau^*_{T,\d})}\left\{\inf_{t \ge
\tau^*_{T,\d}} \(z^+(\tau^*_{T,\d})\;e^{-2\int_{\tau^*_{T,\d}}^t
y^+(s)\, ds}+ \xi
\chi^+_{\tau^*_{T,\d}}(t)\)>\eps\right\}\nonumber\\&&
\hskip10.5cm\le e^{-cT}
\end{eqnarray}
where the last inequality is obtained by the use of Lemma
\ref{x0}.

\vskip.1cm

We prove, now, the statement for the second case
$z^+(\tau^*_{T,\d})<-\eps$, $z^*(\tau^*_{T,\d})>\d$.  At first, we
show that, with large probability, $z^*(t)$ reaches the line
$\frac 3 2 \; t$, i.e. that the stopping time $\tau'_T:=\inf\{t
\ge \tau^*_T: z^*(t)\ge \frac 3 2 \; t\}$ is finite. We compare
$z^*(t)$ with the process $v^+(t)$, solution of the linear problem
\begin{equation}\label{x+-}
d v^+(t)= \frac t 2 \;  v^+(t) \, dt + \xi  dw(t), \quad\quad
v^+(\tau^*_{T,\d})=z^*(\tau^*_{T,\d}),
\end{equation}
we have
\begin{equation}\label{x+-'}
v^+(t)= v^+(\tau^*_{T,\d}) \;e^{\frac 1
4(t^2-\tau_{T,\d}^{*2})}+\xi \; e^{\frac {t^2} 4 }
\int_{\tau^*_{T,\d}}^t e^{-\frac {u^2}4} \; dw_u.
\end{equation}

Since $-z(z+2y^*_t)>tz/2$ for $0\le z \le -3y^*_t/2$,  by Lemma
\ref{lemma:comp}, $z^*(t)\ge v^+(t)$, as long as $0\le z^*(t)\le
\frac 3 2 t$. It is sufficient to apply  the inequality \eqref{MS}
to $v^+(t)$ whose quadratic variation is easily estimable from
\eqref{x+-'} and \eqref{eq:exp-} to show that
\begin{eqnarray*}
&&\mathbf 1_{v^+(\tau^*_{T,\d})>\d}\;\p_{\tau^*_{T,\d},
Y(\tau^*_{T,\d})}\left\{\inf_{t \ge \tau^*_{T,\d}} v^+(t)\le
0\right\}\le e^{-cT} \\ &&\text{and} \quad \mathbf
1_{v^+(\tau^*_{T,\d})>\d}\;\p_{\tau^*_{T,\d},
Y(\tau^*_{T,\d})}\left\{\sup_{t \ge \tau^*_{T,\d}} \frac{2
v^+(t)}{3t}<1\right\}\le e^{-cT}
\end{eqnarray*}
hence
\begin{equation}\label{.1}
\mathbf
1_{\tau^*_{T,\d}<\tau'_T}\;\p_{\tau^*_{T,\d},Y(\tau^*_{T,\d})}\left\{
\tau'_T=\i\right\}=1_{\tau^*_{T,\d}<\tau'_T}\;\p_{\tau^*_{T,\d},Y(\tau^*_{T,\d})}\left\{\sup_{t
\ge \tau^*_{T,\d}} \frac{2 z^*(t)}{3t}<1\right\}\le e^{-cT}
\end{equation}
For $t$ large enough, $y^*(t)\ge t- 1 /t$, thus, from \eqref{.1},
with $\p_{\tau^*_{T,\d},Y(\tau^*_{T,\d})}$-probability greater
than $1-e^{-cT}$,
 there exists $\tau^*_T\le\tau'_T <\i$ such that
$z^+(\tau'_T)\ge -\tau'_T/ 2 -1 /\tau'_T$.

\vskip.1cm

By an analogous comparison argument it is possible to prove that
\begin{equation}\label{z-}
\mathbf 1_{z^+(\tau^*_T)<-\eps, \, \tau'_T<\i}\: \p_{\tau'_T,
z^+(\tau'_T)}\left\{\sup_{t \ge \tau^*_T} z^+(t)<-\eps\right\} \le
e^{-cT}
\end{equation}
\eqref{step3} follows, then, from \eqref{z+}, \eqref{.1} and
\eqref{z-}. \qed

\vskip.6cm

\begin{prop}\label{prop:IV}
There is $c>0$ such that, for any $\eps>0$ small enough, $T,\la$
large enough, $\la<\eps \sqrt T$,
\begin{equation}\label{s3}
\mathbf
1_{\tau^+_{T,\eps}<\i}\p_{\tau^+_{T,\eps},Y(\tau^+_{T,\eps})}
\left\{ \inf_{s \ge \tau^+_{T,\eps}}\sup_{t \ge s}|Y(t)-y^+(t)|\;
\sqrt t>\la \right\}\le e^{-c\la^2}
\end{equation}
\end{prop}

{\bf Proof.} Let us suppose $\tau^+_{T,\eps}<+\i$, thus the
relation \eqref{zz} with $\tau^+_{T,\eps}$ in place of $t_0$ holds,
for $t \ge \tau^+_{T,\eps}$. We apply Lemma \ref{x0} to the
process $\chi^+_{\tau^+_{T,\eps}}(t)$, thus,  by symmetry, we get
\begin{equation}\label{eq:x0'}
\p_{\tau^+_{T,\eps},Y(\tau^+_{T,\eps})}\left\{\sup_{t \ge
\tau^+_{T,\eps}} \: |\chi^+_{\tau^+_{T,\eps}}(t)|\:\sqrt{t}>\frac
\la {2\xi} \right\} \le e^{-c \la^2}
\end{equation}
for any $\la$ large enough. Let us define the stopping time
$\tau''_{T,\eps}:=\inf\{t \ge \tau^+_{T,\eps}: |z^+(t)|>2\eps\}$.
We have
\begin{equation*}
\int_{\tau^+_{T,\eps}}^t e^{-2\int_{u}^t y^+(s)\, ds}\; du \le t^2
e^{-t^2} \int_{\tau^+_{T,\eps}}^t \frac {e^{u^2}}{u^2}\; du \le
\frac c t \wedge 1,
\end{equation*}
thus, with $\p_{\tau^+_{T,\eps},Y(\tau^+_{T,\eps})}$-probability
greater than $1-2e^{-c\la^2}$ we have
\begin{equation}\label{.2}
-\eps e^{-2\int_{\tau^+_{T,\eps}}^t y(s)\, ds} -(c \frac {\eps^2}
t \wedge 1) -\frac{\la} {2\sqrt t} \le z^+(t)\le \eps
e^{-2\int_{\tau^+_{T,\eps}}^t y^+(s)\, ds}+\frac{\la} {2\sqrt t}
\end{equation}
for  $\tau^+_{T,\eps}\le t \le \tau''_{T,\eps}$. Assume $\la <\eps
\sqrt T$, thus, since $\tau^+_{T,\eps}\ge T$,  from \eqref{.2} it
follows that
\begin{equation}\label{.3}
\p_{\tau^+_{T,\eps},Y(\tau^+_{T,\eps})}\{\tau''_{T,\eps} <
+\i\}\le
\p_{\tau^+_{T,\eps},Y(\tau^+_{T,\eps})}\left\{\sup_{\tau^+_{T,\eps}\le
t \le \tau''_{T,\eps}}|z^+(t)|< 2\eps \right\}\le e^{-c \eps^2 T}
\end{equation}
then, by \eqref{.2} and \eqref{.3}, we have
\begin{eqnarray*}
\p_{\tau^+_{T,\eps},Y(\tau^+_{T,\eps})} \left\{ \inf_{s \ge
\tau^+_{T,\eps}}\sup_{t
\ge s}|z^+(t)|\; \sqrt t>\la \right\}\le \p_{\tau^+_{T,\eps},Y(\tau^+_{T,\eps})} \{\tau''_{T,\eps}<+\i\} \\
+\,\p_{\tau^+_{T,\eps},Y(\tau^+_{T,\eps})} \left\{ \inf_{s \ge
\tau^+_{T,\eps}}\sup_{t \ge s}|z^+(t)|\; \sqrt t>\la, \;\Big|\;
\tau''_{T,\eps}=\i \right\}\le 2 e^{-c\la^2}
\end{eqnarray*}
hence \eqref{s3} is proved. \qed

\vskip.6cm

\begin{prop}\label{prop:V}
 There exists $c>0$ such that, for any  $T$ large
enough, $\d,\eps$ small enough,
\begin{equation}\label{step33}
\mathbf 1_{ \Pi>T}\;\mathbf 1_{\tau^*_{T,\d}<\i,
\;z^*(\tau^*_{T,\d})>\d}\;
\p_{\tau^*_{T,\d},Y(\tau^*_{T,\d})}\Big\{Y \notin  E^+ \Big\}\le
e^{-c \eps^2 T}
\end{equation}
\end{prop}

 {\bf Proof.} Assume $\Pi>T$, $\tau^*_{T,\d}<\i$ and $z^*(\tau^*_{T,\d})>\d$ then, for
any $\eps,\d>0$ small enough, $T$ large enough,
\begin{eqnarray*}
\p_{\tau^*_{T,\d},Y(\tau^*_{T,\d})}\Big\{ Y \notin E^+ \Big\}\le
\E\[\mathbf
1_{\tau^+_{T,\eps}<\i}\p_{\tau^+_{T,\eps},Y(\tau^+_{T,\eps})}\left\{
Y \notin E^+ \right\}\]+ \,\p_{\tau^*_{T,\d},Y(\tau^*_{T,\d})}\{
\tau^+_{T,\eps}=+\i \}
\end{eqnarray*}

 thus \eqref{step33} follows from Propositions \ref{prop:III} and \ref{prop:IV}. \qed

\vskip.6cm

{\bf Conclusion of  proof of Theorem $\ref{thm:stoc}$.} Let us
suppose $T>t_0$, thus, from the definition of $\Pi$ and
Proposition \ref{prop:I}, we have
\begin{eqnarray}
\p_{t_0,y_0}\left\{ Y \notin E^+ \cup E^-\right\}= \E\[\mathbf 1_{\Pi>T}\;\p_{T,Y(T)}\left\{Y \notin E^+ \cup E^-\right\}\]\nonumber\\
=\E\[\mathbf 1_{\tau^*_{T,\d}<\i, \; \Pi>
\tau^*_{T,\d}}\p_{\tau^*_{T,\d},Y(\tau^*_{T,\d})}\{Y \notin E^+
\cup E^-\}\]\label{Q0}
\end{eqnarray}
\eqref{Q0} is bounded by
\begin{eqnarray}\label{Q1}
&&\E\[\mathbf 1_{\Pi>T}\;\mathbf 1_{\tau^*_{T,\d}<\i,\;
z^*(\tau^*_{T,\d})<-\d}\p_{\tau^*_{T,\d},Y(\tau^*_{T,\d})}\{ Y \notin E^-\}\]\nonumber\\
&&\hskip.4cm+\E\[\mathbf 1_{\Pi> T}\;\mathbf 1_{\tau^*_{T,\d}<\i,
\;z^*(\tau^*_{T,\d})>\d}\p_{\tau^*_{T,\d},Y(\tau^*_{T,\d})}\{Y
\notin E^+ \}\]\le e^{-c\eps^2 T}
\end{eqnarray}
where the inequality follows from  Propositions \ref{prop:II} and
\ref{prop:V}, and holds  for some $c>0$, for any $T$ large enough,
$\d,\eps$ small enough. The result follows from \eqref{Q1} by
performing the limit for $T\to \i$. \qed

\subsection*{Behavior of $Y(t)$ for $t \to -\i$}

In this part of the Section we will provide some  results for the
behavior of $Y(t)$ for negative $t$, $|t|$ large enough.

\vskip.6cm

\begin{prop}\label{prop:VI}
Let $\p$ be the probability law defined at the beginning of this
Section. There is $c>0$ such that for $T,\la$ large enough,
$\la<\sqrt T$,
\begin{equation}\label{VI}
\p\left\{ \sup_{t \le -T}|Y(t)-y^+(t)|\; \sqrt {|t|}>\la
\right\}\le e^{-c\la^2}
\end{equation}
\end{prop}

{\bf Proof.} $z^+(t)$ satisfies the equation \eqref{zz} even in the
limit as $t_0 \to -\i$. $y^+(t)\to +\i$ and  $z^+(t)\to 0$ for $t
\to -\i$, $\p$-a.s., thus
\begin{equation}\label{.4}
z^+(t)=-\int_{-\i}^t z^{+2}(u)e^{-2\int_{u}^t y^+(s)\, ds}\; du +
\xi \chi^+_{-\i}(t)
\end{equation}
We use Lemma \ref{x0} with $t_0=-\i$ to estimate the behavior of
$\chi^+_{-\i}(t)$, then the proof proceeds specularly to proof of
Proposition \ref{prop:IV}. \qed

\vskip.6cm

\begin{prop}\label{prop:VII}
There is $c>0$ such that for $T,S,\la$ large enough, $S<T$,
$\la<\sqrt S$,
\begin{equation}\label{VII}
\mathbf 1_{|Y(T)-T|<\frac \la  {2\sqrt{T}}}\;\p_{-T,Y(T)}\left\{
\sup_{-T \le t \le -S}|Y(t)-y^+(t)|\; \sqrt {|t|}>\la \right\}\le
e^{-c\la^2}
\end{equation}
\end{prop}

 {\bf Proof.} The proof of \eqref{VII} is almost the
same of Proposition \ref{prop:VI}. \qed

\subsection*{Behavior of $Y(t)$ in bounded intervals}

In this part of the Section we study the behavior of  solutions
$Y(t)$ of \ref{x} starting at time $-T$ from $y:|y-T|\le \eps$,
$\eps$ small enough. We recall that the stopping time $\tau_T \in
[-T,T]$ is the first exit time of $Y(t)$ from the rectangle
$\mathcal R_T$ (see \eqref{tau} and \eqref{eq:xx0}).  Notice that
  the condition $|y-T|\le \eps$ guaranties  $(-T, Y(-T))\in \mathcal R_T$.

\vskip.6cm

\begin{lemma}\label{lemma:X3}
There exists $c>0$ such that, for any $T$ large enough, $\eps$
small enough,
\begin{equation*}
\mathbf 1_{|y-T|\le \eps}\; \p_{-T,y}\Big\{
Y(\tau_{T})=2T\Big\}\le e ^{-cT^2}
\end{equation*}
\end{lemma}

{\bf Proof.} We have $y^+(-T)\ge T$, $y\le T+\eps$, then
$z^+(-T)=y-y^+(-T)\le \eps$, hence, by \eqref{zz},
\begin{equation*}
\mathbf 1_{|y-T|\le \eps}\;\p_{-T,y}\Big\{ z^+(t) \le \eps +\xi
\chi^+_{-T}(t), \: \forall \; t \ge -T\Big\}=1,
\end{equation*}
thus, since $\sup_{-T \le t\le T}y^+(t)< T$, we have
\begin{eqnarray*}
\mathbf 1_{|y-T|\le \eps}\;\p_{-T,y}\Big\{ Y(\tau_{T})=2T\Big\}\le
\mathbf 1_{|y-T|\le \eps}\;\p_{-T,y}\Big\{ \sup_{-T \le t \le
T}z^+(t)\ge T\Big\}\\\le \p_{-T,y}\Big\{ \sup_{-T \le t \le
T}\chi^+_{-T}(t)\ge \frac{T-\eps}\xi\Big\}\le e^{-cT^2}
\end{eqnarray*}

where the last inequality follows from \eqref{MS} and Lemma
\ref{x0}. \qed

\vskip.6cm

\begin{lemma}\label{lemma:X1}
There exists $c>0$ such that, for any $T$ large enough,
\begin{equation}\label{X1}
\p_{\tau_{T},-2T}\left\{
 \Pi \ge \tau_{T}+ T^{-1}
  \right\}\le e^{-cT^3}
\end{equation}

\end{lemma}

{\bf Proof.} Consider the process $\tilde y(t)$,   solution of the
ODE $\eqref{eq:y}$ starting from  $-\frac 3 2 T$ at time
$\tau_{T}$.  Let us consider, now, $\tilde z(t):=Y(t)-\tilde
y(t)$, thus $\tilde z_{\tau_{T}}(\tau_{T})=-\frac T 2$. Using
exactly
 the same arguments used in proof of Proposition \ref{prop:II} to show \eqref{.0}, it
is possible to prove that
\begin{equation}\label{X11}
\p_{\tau_{T},-2T}\left\{\sup_{t \ge \tau_{T}}\tilde z(t)\ge 0
\right\}\le e^{-cT^3}
\end{equation}
$\tilde y(t)$ lies below $y^*(t)$, then, from Proposition
\ref{prop:y}, we know that it explodes to $-\i$. It is easy to
show that $\tilde y(t)$ explodes within $\tau_{T}+T^{-1}$ (see
Lemma 2.3.15 in \cite{Ca}), then \eqref{X1} easily follows from
\eqref{X11}. \qed

\vskip.6cm

{\bf Proof of Proposition \ref{prop:3.1}.} Let $Y(-T)=y$ with
$y:|y-T|\le \eps$, then, from Theorem \ref{thm:stoc} we have
\begin{eqnarray}
\p_{-T,y}\left\{ Y \notin \e^{+}_{T} \cup \e^{-}_{T} \right\}=
\p_{-T,y}\left\{ Y \notin \e^{+}_{T} \cup \e^{-}_{T}, \: Y \in E^+
\cup E^-\right\}\nonumber\\ \le \p_{-T,y}\Big\{Y(T) \in
[-2T,y^+(T)-\eps)\cup (y^+(T)+\eps,2T] , \: Y \in
E^+\Big\}\label{eq:X4'}\\+\,\p_{-T,y}\Big\{Y(\tau_{T})=2T\Big\}\label{eq:X4}
\end{eqnarray}
Lemma $\ref{lemma:X3}$ provides a bound for the probability in
\eqref{eq:X4} that assures its convergence to 0 as $T \to \i$. The
term \eqref{eq:X4'} vanishes as $T \to \i$ since, by Theorem
\ref{thm:stoc}, for any $\eps>0$ small enough,
\begin{equation*}
\p_{-T,y}\Big\{ \inf_{T \ge 0}\sup_{t \ge T} \:\:
|Y(t)-y^+(t)|>\eps \;\Big| \: Y \in E^+\Big\}=1
\end{equation*}
hence \eqref{eq:xx??} follows. From Proposition \ref{prop:IV} and
Lemma \ref{lemma:X1} it follows that
\begin{eqnarray}\label{eq:X5+}
\mathbf 1_{(\tau_{T},Y(\tau_{T}))\in \partial \mathcal
R^\pm_{T}}\; \p_{\tau_{T},Y(\tau_{T})}\left\{Y \notin E^\pm
\right\}\le  e^{-c\eps^2T}
\end{eqnarray}
for some $c>0$, thus  \eqref{eq:xx??+} follows from \eqref{eq:X5+}
and  Theorem \ref{thm:stoc}. \qed

\subsection*{Proof of Theorem \ref{thm:3.2}}

We consider  two processes $Y,\bar Y$ solutions of \eqref{x}
starting from $Y(-T)=y$ and $\bar Y(-T)=\bar y$, with $y,\bar y$
such that $|y -T|\le \eps, |\bar y -T|\le\eps$ for some $\eps>0$
small enough. Without lost of generality, we can suppose $\bar y>
y$. We denote by $\mathcal Q_{-T,y,\bar y}$ the probability law of
the coupled process $(Y(t),\bar Y(t))$ by taking the same noise
for $Y$ and $\bar Y$.

Let us fix $S\in (1,T)$ and  $\eps>0$ small enough and define the
sets
\begin{eqnarray}\label{A}
\mathcal A_{S}=\mathcal A_{T,S,\eps}:=\Bigg\{Y\; : \: \sup_{-T\le
t \le -S} |Y(t)-y^+(t))|\le \eps\Bigg\}
\end{eqnarray}
and
\begin{eqnarray}\label{B}
\mathcal B_{S}=\mathcal B_{S,\eps}:=\Bigg\{Y\; : \: \sup_{t \ge S}
|Y(t)-y^+(t))|\le \eps\Bigg\}
\end{eqnarray}
Let $\tau_S, \bar \tau_S$ be the first exit times respectively for
the processes $Y(t)$ and $\bar Y(t)$ from the rectangle $\mathcal
R_S$ defined in \eqref{eq:xx0}. We call $\Pi$ and $\bar \Pi$ the
times of explosion to $-\i$ of $Y$ and $\bar Y$.  $Y(t)$ and $\bar
Y(t)$ are well defined, thus, respectively for $t \le \Pi$ and for
$t \le \bar \Pi$. We agree with the convention to define
$Y(t):=-\i$ for $t \ge \Pi$, $\bar Y(t):=-\i$ for $t \ge \bar
Y(t)$. We have the following results.

\vskip.6cm

\begin{lemma}\label{lemma:Y3}
For any $S\in (1,T)$ and  $\eps>0$ small enough
\begin{equation}\label{Q2}
Q_{-T,y,\bar y}\left\{\lim_{T \to \i}\sup_{t \ge -T}\bar Y(t)-Y(t)=0\; \Big|
\: Y,\bar Y \in \mathcal A_{S} \cap\mathcal E_S^+ \cap \mathcal
B_S \right\}=1.
\end{equation}
\end{lemma}

 {\bf Proof.}
Let us assume  $Y,\bar Y \in \mathcal A_S  \cap\mathcal E_S^+ \cap
\mathcal B_S $. We denote by $v(t)$ the process $\bar Y(t)- Y(t)$,
then $dv=-v(Y+\bar Y)dt$, hence
\begin{equation}\label{x'}
v(t)=(\bar y- y)e^{-\int_{-T}^t(Y(u)+\bar Y(u))du}
\end{equation}
 thus $v(t)>0$ for any $t \ge -T$.\\
 Since $Y,\bar Y \in \mathcal A_S$ and $y^+(t)\ge -t$
for $t<0$, from \eqref{x'} we have
\begin{equation*}
0 \le v(t)\le  \eps \: e^{2\int_{-T}^t (u+ \eps)\,du}\le \eps
\,e^{(S-\eps)^2}\,e^{-(T-\eps)^2} \quad \text{for}\quad -T \le t
\le -S
\end{equation*}
thus
\begin{equation}\label{.8'}
\lim_{T \to \i}\sup_{-T \le t \le -S}v(t)=0.
\end{equation}

$Y,\bar Y \in \mathcal E^+_S$ implies $\bar Y(t),Y(t)\ge -2S$ for
any $-S \le t \le S$, then, by \eqref{x'},
\begin{equation*}
0\le v(t)\le v(-S)e^{-\int_{-S}^t(Y(u)+\bar Y(u))du}\le v(-S)
e^{8S^2} \quad \text{for} \quad -S\le t \le S,
\end{equation*}
thus, by \eqref{.8'},
\begin{equation}\label{.8''}
\lim_{T \to \i}\sup_{-S \le t \le S}v(t)=0.
\end{equation}

$Y,\bar Y \in \mathcal B_S$ thus $\bar Y(t), Y(t)\ge
y^+(t)-\eps>0$, for $t\ge S$, $S$ large enough, then, from
\eqref{x'} we have
\begin{equation*}
0\le v(t)\le v(S) e^{-\int_{S}^t(Y(u)+\bar Y(u))du}\le v(S) \quad
\text{for}\quad t \ge S
\end{equation*}
thus, from \eqref{.8''},
\begin{equation}\label{.8'''}
\lim_{T \to \i}\sup_{t \ge S}v(t)=0
\end{equation}
then the Lemma is proved. \qed

\vskip.6cm

\begin{lemma}\label{lemma:Y4}
For any $S\in (1,T)$ and $\eps>0$ small enough
\begin{equation}\label{Q3}
Q_{-T,y,\bar y}\left\{\lim_{T \to \i}\sup_{-T \le t \le
\tau_S}|\bar Y(t)-Y(t)|=0\; \Big| \: Y \in \mathcal A_{S}\cap
\mathcal E_S^-, \; \bar Y \in \mathcal A_S \right\}=1.
\end{equation}
\end{lemma}

 {\bf Proof.}
Let us assume $Y \in \mathcal A_{S}\cap \mathcal E_S^-$ and $\bar
Y \in \mathcal A_S $. Consider the process $v(t)$ defined in the
proof of the previous Lemma, then, since $Y, \bar Y \in \mathcal
A_S$, \eqref{.8'} holds also in the current case.\\
 On the other
hand $Y \in \mathcal E^-_S$ implies $\bar Y(t)\ge Y(t)\ge -2S$,
then, since $|\tau_S| \le S$,
\begin{equation*}
0\le v(t)\le v(-S) e^{-\int_{-S}^t(Y(u)+\bar Y(u))du}\le
v(-S)e^{8S^2} \quad \text{for} \quad -S\le t \le \tau_S,
\end{equation*}
thus, by \eqref{.8'}, we have
\begin{equation}\label{Q4}
\lim_{T \to \i}\sup_{-S \le t \le \tau_S}v(t)=0  \quad \quad
Q_{-T,y,\bar y}-\text{a.s.}
\end{equation}
hence \eqref{Q3} follows. \qed

\vskip.6cm

\begin{prop}\label{lemma:Y1}

For any bounded continuous function $g(y)$ with compact support
and for any fixed $\eps>0$ small enough, $t \ge -T$, we have
\begin{equation}\label{eq:Y1}
\lim_{T \to \i}\mathbf 1_{|y-T|\le \eps}\mathbf 1_{|\bar y-T|\le
\eps }\bigg|\E_{\p_{-T,y}}\[g(Y(t))\]-\E_{\p_{-T,\bar y}}\[g(\bar
Y(t))\]\bigg|=0.
\end{equation}
\end{prop}

 {\bf Proof.}
We define $G(t):=\big|g(Y(t))-g(\bar Y(t))\big|$ then we need to
prove that
\begin{equation}\label{eq:E0}
\lim_{T \to \i}\mathbf 1_{|y-T|\le \eps} \mathbf 1_{|\bar y-T|\le
\eps}\E_{\mathcal Q_{-T,y,\bar y}}\[G(t)\]=0.
\end{equation}

Let us fix $S \in (1,T)$ large enough and $\eps>0$ small enough,
$y,\bar y:|y-T|\le \eps, |\bar y-T|\le \eps$. For $\mathcal A_S$
  as in \eqref{A} we have
\begin{eqnarray}\label{.10}
\big|\E_{\mathcal Q_{-T,y,\bar y}}\[G(t)\]-\E_{\mathcal
Q_{-T,y,\bar y}}\[\mathbf 1_{Y,\bar Y \in \mathcal A}G(t)\]\big|
\hskip2cm\\\le
 2\sup |g| \(\p_{-T,y}\left\{Y \notin \mathcal A_S\right\}+\p_{-T,\bar y}\left\{\bar Y \notin \mathcal A_S\right\}\)\nonumber\\
 \le 4  \sup |g| \;e^{-c \eps^2 S}\nonumber
\end{eqnarray}
where the last inequality follows from \eqref{VII}. We have
\begin{eqnarray}\label{.12}
&&\Big|\E_{\mathcal Q_{-T,y,\bar y}}\[\mathbf 1_{Y,\bar Y \in
\mathcal A_S}G(t)\]- \E_{\mathcal Q_{-T,y,\bar y}}\[\mathbf 1_{Y,
\bar Y \in\mathcal A_S \cap (\mathcal E^+_S \cup \mathcal E^-_S)}
G(t)\] \Big| \\&&\le 2\sup |g|\(\p_{-S,y}\{Y \notin \mathcal E^+_S
\cup \mathcal E^-_S\}+ \p_{-S,\bar y}\{\bar Y \notin \mathcal
E^+_S \cup \mathcal E^-_S\}\)\nonumber
\end{eqnarray}

For $\mathcal B_S$ as in \eqref{B} we have
\begin{eqnarray}\label{.13}
\Big|\E_{\mathcal Q_{-T,y,\bar y}}\[\mathbf 1_{Y,\bar Y \in
\mathcal A \cap \mathcal E^+_S}G(t)\]-\E_{\mathcal Q_{-T,y,\bar
y}}\[\mathbf 1_{Y,\bar Y \in \mathcal A \cap \mathcal B_S \cap
\mathcal E^+_S}G(t)\]\Big|\hskip1.3cm\\ \le 2\sup|g| \(\mathbf
1_{|y-y^+(S)|\le \eps}\p_{S,y}\{Y \notin \mathcal B_S\}+ \mathbf
1_{|\bar y-y^+(S)|\le \eps}\p_{S,\bar y}\{\bar Y \notin \mathcal
B_S\}\)\nonumber
\\\le 4 \sup |g| \;e^{-c\eps^2 S}\nonumber
\end{eqnarray}
where the last inequality follows from \eqref{s3}.

Since $g$ is bounded, it follows from  \eqref{.10}, \eqref{.12},
\ref{eq:xx??} and \eqref{.13} that for any $\zeta>0$ there exists
$S_0$ such that, for any $T>S\ge S_0$, for any $t \ge -T$,
\begin{equation}\label{Q5}
\Big|\E_{\mathcal Q_{-T,y,\bar y}}\[G(t)\]- \E_{\mathcal
Q_{-T,y,\bar y}}\[\(\mathbf 1_{Y, \bar Y \in\mathcal A_S \cap
\mathcal E^+_S \cap \mathcal B_S} + \mathbf 1_{Y, \bar Y
\in\mathcal A_S \cap \mathcal E^-_S}\)G(t)\]\Big|\le \zeta.
\end{equation}

By the continuity of $g$ and Lemma \ref{lemma:Y3} it follows that
\begin{equation}\label{Q2'}
Q_{-T,y,\bar y}\left\{\lim_{T \to \i}\sup_{t \ge -T}G(t)=0\; \Big|
\: Y,\bar Y \in \mathcal A_{S} \cap\mathcal E_S^+ \cap \mathcal
B_S \right\}=1,
\end{equation}
thus
\begin{equation}\label{.14''}
\lim_{T \to \i}\mathbf 1_{t \ge -T}\E_{\mathcal Q_{-T,y,\bar
y}}\[\mathbf 1_{Y,\bar Y \in \mathcal A_S \cap \mathcal B_S \cap
\mathcal E^+_S}G(t)\]=0.
\end{equation}
On the other hand, By Lemma \ref{lemma:Y4},
\begin{equation}\label{Q33}
Q_{-T,y,\bar y}\left\{\lim_{T \to \i}\sup_{-T \le t \le
\tau_S}G(t)=0\; \Big| \: Y \in \mathcal A_{S}\cap \mathcal E_S^-
\; \bar Y \in \mathcal A_S \right\}=1,
\end{equation}
 thus
\begin{equation}\label{.14}
\lim_{T \to \i}\E_{\mathcal Q_{-T,y,\bar y}}\[\mathbf 1_{-T\le t
\le \tau_S}\mathbf 1_{Y,\bar Y \in \mathcal A_S  \cap \mathcal
E^-_S}G(t)\]=0.
\end{equation}

We have
\begin{equation}\label{.15}
\E_{\mathcal Q_{-T,y,\bar y}}\[\mathbf 1_{\tau_S \le t \le \tau_S
+S^{-1}}\mathbf 1_{Y,\bar Y \in \mathcal A \cap \mathcal
E^-_S}G(t)\]\le 2 \sup |g|\; \p_{-T,y}\{t-S^{-1}\le \tau_S \le t\}
\end{equation}
with the right hand side term  vanishing as $S \to \i$. \\
Since $\bar Y(t)\ge Y(t)$, $\bar \Pi \ge \Pi$, thus, for $t \ge
\bar \Pi$, $Y(t)=\bar Y(t)=-\i$, then $G(t)=0$. It remains to
estimate the term for $\tau_S +S^{-1}\le t \le \bar \Pi$. We have
$Y(\tau_S)=-2S$ and $\bar Y(\tau_S)=-2S+ v(\tau_S)$, with, by
Lemma \ref{lemma:Y4}, $\lim_{T \to \i}|v(\tau_S)|=0$ for $Y,\bar Y
\in \mathcal A_S \cap \mathcal E^-_S$. Hence for any fixed
$\zeta>0$ arbitrarily small there is $T_0$ such that, for any
$T>T_0$
\begin{eqnarray*}
&&\E_{\mathcal Q_{-T,y,\bar y}}\[\mathbf 1_{\tau_S +S^{-1}\le t
\le \bar \Pi}\;\mathbf 1_{Y,\bar Y \in  \mathcal A_S \cap
\mathcal E^-_S}G(t)\]\\
&&\le \E_{\mathcal Q_{-T,y,\bar y}}\[\mathbf 1_{\tau_S +S^{-1}\le
\bar \Pi}\;\mathbf 1_{\bar Y(\tau_S)\le -2S+\zeta}\; G(t)\]
\\&&\le 2\sup |g| \; \E_{\mathcal Q_{-T,y,\bar
y}}\[\mathbf 1_{\bar Y(\tau_S)\le -2S+\zeta}\;\p_{\tau_S,\bar
Y(\tau_S)}\left\{\bar \Pi \ge \tau_S +S^{-1}\right\}
\]
\end{eqnarray*}
then, by Lemma \ref{lemma:X1}, for any $\zeta'>0$ there exists
$S_0$ such that, for any $T>S>S_0$,
\begin{eqnarray}\label{.16'}
\E_{\mathcal Q_{-T,y,\bar y}}\[\mathbf 1_{\tau_S +S^{-1}\le t \le
\bar \Pi}\;\mathbf 1_{Y,\bar Y \in  \mathcal A_S \cap \mathcal
E^-_S}G(t)\] < \zeta'.
\end{eqnarray}
From \eqref{.15} and \eqref{.16'} it follows that
\begin{equation}\label{.14'}
\lim_{T \to \i}\E_{\mathcal Q_{-T,y,\bar y}}\[\mathbf 1_{t\ge
\tau_S}\mathbf 1_{Y,\bar Y \in \mathcal A_S  \cap \mathcal
E^-_S}G(t)\]=0,
\end{equation}

then \eqref{eq:E0} follows from \eqref{Q5}, \eqref{.14''},
\eqref{.14} and \eqref{.14'}. \qed

\vskip.6cm

\begin{coro}\label{coro:Y2}
Let $Y(t),\hat Y(t)$ be solutions of \eqref{x} starting from
$Y(-T)=y$, $\hat Y(-S)= \hat y$, $T>S$, then, for any function
$g(y)$ as in the previous Proposition,
\begin{equation}\label{eq:Y2}
\lim_{S \to \i}\lim_{T \to \i}\Big|\E_{\p_{-T,y}}[g(Y(t))]\mathbf
1_{|y-T|\le \eps}-\E_{\p_{-S,\hat y}}[g(\hat Y(t))]\mathbf
1_{|\hat y-S|\le \eps}\Big|=0
\end{equation}
\end{coro}

 {\bf Proof.} Suppose $T>S$, $|y-T|\le \eps$, $|\hat y -S|\le
 \eps$. We have
\begin{equation}\label{Q6'}
\E_{\p_{-T,y}}\[g(Y(t))\]=\E_{\p_{-T,
y}}\[\E_{\p_{S,Y(-S)}}\[g(Y(t))\]\]
\end{equation}
thus
\begin{eqnarray*}
\bigg|\E_{\p_{-T,y}}\[g(Y(t))\]-\E_{\p_{-T,y}}\[\E_{\p_{-S,Y(-S)}}\[g(Y(t))\]\mathbf
1_{|Y(-S)-S|\le \eps}\]\bigg| \\\le \sup |g| \;
\p_{-T,y}\left\{|Y(-S)-S|>\eps\right\}\le \sup|g|\;
e^{-c\eps^2S}\label{Q6}
\end{eqnarray*}
where the last inequality in \eqref{Q6} follows from \eqref{VII}.
On the other hand, from Proposition \ref{lemma:Y1}, for any
$\zeta>0$ there exists $S_0$ such that, for any $T>S>S_0$,
\begin{equation}\label{Q6''}
\bigg|\E_{\p_{-T,y}}\[\E_{\p_{S,Y(-S)}}\[g(Y(t))\]\mathbf
1_{|Y(-S)-S|\le \eps}\]-\E_{\p_{-S,\hat y}}\[g(\hat Y(t))\]\mathbf
1_{|\hat y-S|\le \eps}\bigg|<\zeta
\end{equation}
then \eqref{eq:Y2} follows from \eqref{Q6}, \eqref{Q6'},
\eqref{Q6''} and the boundedness of $g$. \qed

\vskip.6cm

\begin{prop}\label{prop:IX}
Let $\p$ be the probability law defined at the beginning of this
Section, then the probabilities $p_\pm :=\p\{Y \in E^\pm\}$ are
strictly positive.
\end{prop}

{\bf Proof.} Let us prove, at first, the statement for $E^-$. By
\eqref{.4}, $z^+(t)\le \xi \chi^+_{-\i}(t)$  $\p$-a.s., thus, for
$\gamma:=y^+(0)-y^*(0)>0$, $\Phi(x):=\frac{1}{\sqrt {2
\pi}}\int_{x}^{+\i} e^{-\frac{z^2}2}\, dz$, we have
\begin{eqnarray}\label{.5}
\p\{z^*(0)<0\}= \p\{z^+(0)<-\gamma\}\ge \p\{\chi^+_{-\i}(0)<-
\gamma \xi^{-1}\}=\Phi \({\gamma}\xi^{-1}\Big/
\sqrt{\E\[\chi^{+2}_{-\i}(0)\]}\)\ge c >0\nonumber
\end{eqnarray}
since, by Lemma \ref{x0}, $\E\[\chi^{+2}_{-\i}(0)\]$ is bounded by
a constant. From \eqref{eq:z'} it is easy to verify that
\begin{equation*}
z^*(t)\le z^*(0)\,e^{-2\int_{0}^t y^*(s)\, ds} + \xi
\chi_{0}^*(t),  \quad \quad \chi^*_{0}(t):=\int_{0}^t
e^{-2\int_u^t y^*(s)\, ds}\; dw_u
\end{equation*}
 $\p$-a.s. Let us suppose  $z^*(0)<0$, then  $z^*(t)\le \xi \chi_{0}^*(t)$,
thus
\begin{eqnarray}\label{.6}
\p\{z^*(T)<-\d\}\ge \p\{ \chi_{0}^*(T)>- \d \xi^{-1}\}= \Phi\(\d
\xi^{-1}\sqrt{\E\[\chi^{*2}_{0}(T)\]}\)\ge c
>0.\nonumber
\end{eqnarray}
for any $\d>0$, hence, from Proposition \ref{prop:II}, \eqref{.5}
and \eqref{.6} it follows that $\p\{Y \in E^-\}>0$.

\vskip.1cm

By the use of comparison arguments as in proof of Proposition
\ref{prop:IV} it is easily provable that there exist
$\eps_{\text{max}},c>0$ such that, for any $\eps \le
\eps_{\text{max}}$, if $|z^+(-T)|<\frac \eps 2$, then
\begin{equation}\label{.7}
\p_{-T,Y(-T)}\left\{\sup_{-T\le t \le T}|z^+(t)|<\eps\right\}>
e^{-\frac c {\eps^2}}
\end{equation}
To prove \eqref{.7} it is sufficient to use the small balls
inequality \eqref{SB}. Thus the claim for $E^+$ follows from
Proposition \ref{prop:IV}, Proposition \ref{prop:VI} and
\eqref{.7}. \qed

\vskip.6cm

{\bf Conclusion of the proof of Theorem \ref{thm:3.2}.} The
convergence result \eqref{eq:Y?} is a direct consequence of
Corollary \ref{coro:Y2}.  \eqref{eq:Y??} easily follows from
Proposition \ref{prop:VI}. \\The convergence of the probabilities
$\mathbf 1_{|y-T|\le \eps}\p_{-T,y}\{Y \in E^\pm\}$ is a direct
consequence of \eqref{eq:Y?}; finally,  from \eqref{eq:xx??+} and
Proposition \ref{prop:IX} it follows that $p_-\in(0,1)$. \qed

\section{Escape from criticality} \label{section:escape}

In this Section we study our $N$-finite dynamics assuming $Y_N \in
\e_{T}^-$, i.e. $Y_N(\tau_{N,T})=-2T$ or, equivalently,
\begin{equation}\label{hp}
m_N(\mathcal T_{N,T})= x_c - \frac{2T} {\nu N^{1/ 3}}, \quad\quad
\mathcal T_{N,T}:=\mu \, \tau_{N,T}N^{1/3}\in [- T, T]\mu N^{1/3}
\end{equation}
we recall the definition of $\mathcal H^\pm_\ga(I)$ in \eqref{H+-}
and  prove the following result.

\vskip.5cm

\begin{prop}\label{prop:escape}
For any  $\mu'>\mu$
\begin{equation}\label{eq:escape}
\lim_{T \to \i}\lim_{N \to \i}\p_{N}\left\{ \mathcal H^-_0\(\{\mu'
T N^{1/3}\}\)\,\bigg| \: Y_N \in \e_{T}^-\right\}=1
\end{equation}
\end{prop}

 We consider the stochastic process
 $x_N^*(t):=x_N^*(\mathcal T_{N,T})(t)$ defined as the solution of the ODE
$\eqref{eq:eqF}$ with $h=h_N(t)$ and random  initial condition
$x_N^*(\T_{N,T})=m_N(\T_{N,T})$. We prove that, for any $\mu'>\mu$
($\mu$ as in \eqref{mu}), $x_N^*(t)$ reaches $X_-(h_N(t))$ within
the time $\T_{N,T}+\mu' T N^{1/3}$, then we show that, by tracking
$x_N^*(t)$, our magnetization $m_N(t)$ approaches $X_-(h_N(t))$.
We denote by $\p^-_{N,T}$ the probability law of $m_N(t)$ given
$Y_N \in\e_{T}^-$.  All the computations are done for $N>T$, $N,T$
large enough. Unless further indications, we will denote by $c$ a
generic positive constant independent of $N,T$. In order to
lighten notation, in this Section we will omit the index  $N$ for
the magnetization and simply write $m(t)$ and $x^*(t)$.

 \vskip.3cm

We define  the stopping time
\begin{equation*}
\hat {\mathcal T}_{N,T}:=\inf \left\{t \ge\mathcal T_{N,T} :
|m(t)-x^*(t)|> N^{-1/6}\right\}
\end{equation*}
 and recall that $F(x)=-x+ \tanh\{\b(x+h)\}$, we have the
following result

 \vskip.6cm

\begin{lemma}\label{lemma:esc}
Let $\tau,\tau'$ be two stopping times for $m(t)$ such that
 $ \T_{N,T}<\tau<\tau'<\hat \T_{N,T}$ and $N>\tau'-\tau$  $\p^-_{N,T}$-a.s.  There exists a
function $\psi(t)$, such that
\begin{eqnarray}\label{eq:B2(4)'}
\sup_{\T_{N,T}<t<\hat \T_{N,T}}\Big|\psi(t)-\frac {\partial
F}{\partial x}(x^*(t),h_N(t))\Big|\le  c N^{-1/6}
\end{eqnarray}
and, for   $\gamma>0$ small enough,
\begin{equation}\label{eq:B2(6)} \p^-_{N,T}\bigg\{\sup_{ \tau \le t \le \tau'}\[|m(t)-x^*(t)|-
\Theta_{\tau,\tau'}(t)\]\le 0 \bigg\}\ge1-cN^{-\gamma}
\end{equation}
with
\begin{equation}\label{Theta}
 \Theta_{\tau,\tau'}(t):=|m(\tau)-x^*(\tau)|e^{\int_{\tau}^t \psi(u)du}+
 \frac{2(\tau'-\tau)^{1/ 2}}{N^{(1 - \gamma)/ 2}}\(1+e^{\int_{\tau}^t
\psi(u)du}\int_{\tau}^t |\psi(s)|e^{-\int_{\tau}^s \psi(u)du}\,
ds\)\nonumber
\end{equation}

\end{lemma}

{\bf Proof.} Let us define the function $f(x,t):=x-x^*(t)$, then
the process
\begin{equation*}
\mathcal
M(t):=f(m(t),t)-f(m(\T_{N,T}),\T_{N,T})-\int_{\T_{N,T}}^{t}\[\L_{h(s)}f+\frac{\partial
f}{\partial s}\](m(s),s)\, ds
\end{equation*}
is a martingale. For any $\tau$ as in the hypothesis,  the process
$\M_\tau(t):=\M(t)-\M(t \wedge \tau)$ is a martingale as well and
\begin{equation*}
\mathcal V_\tau(t):=\int_{t \wedge \tau}^t
[\L_{h(s)}f^2-2f\L_{h(s)} f](m(s),s)\, ds
\end{equation*}
is its quadratic variation. For any $t\ge\T_{N,T}$,
\begin{eqnarray*}
 \E\[\M_{\tau}^2(t \wedge \tau')|\,\tau,\tau'\]=\E\[\mathcal V_{\tau}(t \wedge \tau')|\,\tau,\tau'\]
 \le c N^{-1} \(t \wedge \tau'-t \wedge \tau\) \le
 c  N^{-1}(\tau'-\tau)
\end{eqnarray*}
thus, by the Doob's inequality, for any $\gamma \in (0,1)$,
\begin{eqnarray*}
&&\p\left\{ \sup_{\tau \le t \le \tau'}|\M_{\tau}(t)|\ge
\frac{(\tau'-\tau)^{1/ 2}}{N^{(1 - \gamma)/ 2}}\; \bigg|\:
\tau,\tau'\right\}\le c N^{-\gamma}
\end{eqnarray*}
then
\begin{eqnarray}
&&\hskip-1cm\p^-_{N,T}\bigg\{\sup_{\tau \le t \le \tau'}|\M_{\tau}(t)|\ge \frac{(\tau'-\tau)^{1/ 2}}{N^{(1 - \gamma)/ 2}}\bigg\}\label{eq:B2(1)} \\
&&=\E_{\p^-_{N,T}}\[\mathbf 1_{\tau<\tau'}\;\p\left\{ \sup_{\tau
\le t \le \tau'}|\M_{\tau}(t)|\ge \frac{(\tau'-\tau)^{1/ 2}}{N^{(1
- \gamma)/ 2}}\; \bigg|\: \tau,\tau'\right\}\]\le c
N^{-\gamma}\nonumber
\end{eqnarray}

Recall the initial condition \eqref{hp},  then there exists a
function  $\psi(t)$ satisfying  \eqref{eq:B2(4)'} and such that
\begin{eqnarray}\label{eq:B2(4)}
\Big|\big[\L_{h(s)}f+\frac{\partial f}{\partial
s}\big](m(s),s)-\psi(s)(m(s)-x^*(s))\Big|\le cN^{-1}
\end{eqnarray}
for $t \ge \T_{N,T}$. For $\tau\le t \le \hat\T_{N,T}$ we define
the process
\begin{equation}\label{Q22}
R_{\tau}(t):=f(m(t),t)- f(m(\tau),\tau)-\int_{\tau}^t
\psi(s)f(m(s),s)\,ds -\M_{\tau}(t),
\end{equation}
then, from \eqref{eq:B2(4)},
\begin{equation}\label{Q21}
\sup_{ \tau \le  t \le \hat
\T_{N,T}}\frac{|R_{\tau}(t)|}{|t-\tau|}\le c\, N^{-1}\quad \quad
\p^-_{N,T}-\text{a.s.}
\end{equation}
 By treating $\eqref{Q22}$ as an  integral equation
for $f(m(t),t)$ we find
\begin{eqnarray*}
&&\hskip-1.5cm f(m(t),t)=f(m(\tau),\tau)\,e^{\int_{\tau}^t \psi(u)du} +[R_{\tau}(t)+\M_{\tau}(t)] \\
&&\hskip.3cm  +\,e^{\int_{\tau}^t \psi(u)du}\int_{\tau}^t
[R_{\tau}(s)+\M_{\tau}(s)]\,\psi(s)\,e^{-\int_{\tau}^s
\psi(u)du}\, ds, \quad \quad \tau \le t \le \hat\T_{N,T}
\end{eqnarray*}
From  $\eqref{eq:B2(1)}$ and \eqref{Q21}, assuming $N>\tau'-\tau$,
we find
\begin{equation}\label{eq:B2(0)}
\p^-_{N,T}\bigg\{\sup_{\tau\le t \le
\tau'}|R_{\tau}(t)+\M_{\tau}(t)|\le   \frac{2(\tau'-\tau)^{1/
2}}{N^{(1 - \gamma)/ 2}}\bigg\} \ge 1-cN^{- \gamma}
\end{equation}
thus \eqref{eq:B2(6)} follows. \qed

\vskip.6cm

\begin{lemma}\label{lemma:Q1}

Let us fix $\d>0$ small enough and consider the stopping time
 \begin{equation*}
  \T_{N,T}'=\T_{N,T,\d}':= \inf\big\{t \ge \T_{N,T}: x^*(t)\le x_c-\d\big\}
\end{equation*}
then, $\p^-_{N,T}$-a.s, there exists $C_0>0$ such that
$\T_{N,T}'-\T_{N,T}\le C_0 \,T^{-1}N^{\frac 1 3}$ for any $T,N$
large enough, and
\begin{equation}\label{esc2}
 x^*(t)\ge  \hat x(t):=x_c -\frac {2T} {\nu N^{\frac 1 3}- 4\b x_c T(t+T\mu N^{\frac 1
   3})} \quad \quad\text{for}\quad \T_{N,T}\le t \le\T_{N,T}'
\end{equation}

\end{lemma}

{\bf Proof.} Recall the initial condition \eqref{hp} for $x^*$,
then, for $ \T_{N,T} \le t \le 2\mu T N^{\frac 1 3}$ we have $-\hc
\le h_N(t)\le -\hc+ c T^2N^{-\frac 2 3}$,
 thus
\begin{eqnarray*}
0\le F(x^*(t),h_N(t))-F(x^*(t),h_c)\le cT^2N^{-\frac 2 3}
\end{eqnarray*}
 hence, for $\T_{N,T} \le t \le  \T_{N,T}' \wedge 2\mu TN^{\frac 1 3}$,
\begin{eqnarray*}
-\b x_c(1+c_0\d)(x^*(t)-x_c)^2 \le F(x^*(t),h_N(t))\le -\b
x_c(1-c_0\d)(x^*(t)-x_c)^2 +cT^2N^{-\frac 2 3}
\end{eqnarray*}

for a suitable $c_0>0$ independent of $\d$, then  $x_1(t)\le
x^*(t)\le x_2(t)$,  $x_{1,2}(t)$ solutions of
\begin{eqnarray}\label{eq:m1,2}
&&\hskip-1cm x_1'(t)= -\b x_c (1+c_0\d)(x_1(t)-x_c)^2, \\
&& x_2'(t)= -\b x_c(1-c_0\d)(x_2(t)-x_c)^2+cT^2N^{-\frac 2 3}
\end{eqnarray}

with $x_1(\T_{N,T})=x_2(\T_{N,T})=x^*(\T_{N,T})=x_c-2T/\nu N^{1/
3}$. It is easy to check that
 \begin{equation}\label{eq:m1'}
   x_1(t)=x_c -\frac 1 {\nu N^{1/3}(2T)^{-1}- \b x_c (1+c_0
   \d)(t-\tau_0)}.
\end{equation}
On the other hand  $m_2(t)$ is a function blowing up at time
\begin{equation*}
\T_{N,T}+ C_0\; \frac {N^{1/ 3}}T<< 2\mu T N^{1/ 3}
\end{equation*}
for a suitable $C_0$ possibly depending on $\d$. In particular we
have $\T_{N,T}' < 2\mu TN^{1 /3} $, thus the result follows. \qed

\vskip.6cm

\begin{lemma}\label{lemma:Q2}
Let us fix $\d>0$ small enough and define  the stopping time
 \begin{equation*}
 \T_{N,T}''= \T_{N,T,\d}'':=\inf\big\{t \ge \T_{N,T}': x^*(t)\le X_-(h_N(t))+\d\big\}
\end{equation*}
then,  $\p^-_{N,T}$-a.s., there exists $C_1>0$ such that
$\T_{N,T}'' - \T_{N,T}'\le C_1$ for any $T,N$ large enough, and
 \begin{equation}\label{stepII}
x^*(t)\le X_-(h_N(t))+\d \quad\quad  \text{for any}\quad t \ge
\T_{N,T}''
\end{equation}

\end{lemma}

{\bf Proof.}  Let $\bar x_1(t)$ and $\bar x_2(t)$ be the solutions
of
\begin{equation}\label{eq:m1,2'}
\bar x'_1(t)= F(\bar x_1(t),h_c) \quad \text{and} \quad \bar
x'_2(t)= F(\bar x_2(t),h_c)+cT^2N^{-\frac 2 3}
\end{equation}
with $\bar x_1(\T_{N,T}')=\bar x_2(\T_{N,T}')=x^*(\T_{N,T}')= x_c
-\d$.  From Lemma \ref{lemma:Q1} we know that  $\T_{N,T}\le 2\mu
TN^{\frac 1 3}$, $\p^-_{N,T}$-a.s., then  $\bar x_1(t)\le
X_+(h_N(t))\le
\bar x_2(t)$ for $\T_{N,T}'\le t \le \T_{N,T}'+\mu TN^{1/ 3}$.\\
Consider the stopping time $\tilde\T_{N,T}'':=\inf\big\{t \ge
\T_{N,T}' :\bar x_1(t)\le X_-(h_N(t))+\d/2\big\}$, then there exists
$C_1>0$ such that $\tilde\T_{N,T}''-\T_{N,T}'\le C_1$. We denote
by $\Delta \bar x(t)$ the nonnegative function $ \bar x_2(t)-\bar
x_1(t)$, thus
\begin{equation*}
\frac d {dt}\Delta \bar x(t)\le (\b-1)\Delta \bar x(t)+c \,
\frac{T^2}{N^{2 /3}} \quad  \quad \Delta \bar x(\T_{N,T}')=0
\end{equation*}

hence  $\Delta \bar x(t)\le c T^2 N^{-2/ 3}$ for any $t \le
\T_{N,T}'+\mu T N^{1/3}$, then, in particular,
$x^*(\tilde\T_{N,T}'')\le \bar x_2(\tilde\T_{N,T}'')\le
\bar x_1(\tilde\T_{N,T}'')+cT^2 N^{- 2/3}\le X_-
(h_N(\tilde\T_{N,T}''))+\d$ for $N$ large enough, then $\T_{N,T}''
\le \tilde\T_{N,T}'' \le \T_{N,T}'+C_1$. \eqref{stepII} is thus
proved. \qed

\vskip.6cm

\begin{lemma}\label{lemma:Q3}
Consider the stopping time
\begin{equation*}
\T_{N,T}''':=\inf\left\{t \ge \T_{N,T}'' : x^*(t)\le
X_-(h_N(t))+N^{-1/2}\right\}
\end{equation*}
then, $\p^-_{N,T}$-a.s., there exists $C_2>0$ such that
$\T_{N,T}''' -\T_{N,T}''\le C_2\ln N$ for any $T,N$ large enough,
and
\begin{eqnarray}\label{stepIII}
|x^*(t)- X_-(h_N(t))|\le N^{-1 /2}  \quad \quad \text{for
any}\quad \T_{N,T}'''\le t \le \frac \pi 2\;  N^{2 /3}
\end{eqnarray}
\end{lemma}

{\bf Proof.}
 There exists  $c>0$ such that
\begin{equation*}
-\frac {\partial}{\partial x}F(X_-(h_N(t)),h_N(t))=\b
[X_-(h_N(t))^2-x_c^2]\ge c
\end{equation*}
for $\T_{N,T}'' \le t \le \frac \pi 2 N^{2/3}$. We have $\d\ge
x^*(t)-X_-(h_N(t))\ge 0$ for $t \ge \T_{N,T}''$, then there exists
$c_0>0$ not depending on $\d$ such that
\begin{equation*}
\F_N(x^*(t),h_N(t))\le -c(1-c_0\d)(x^*(t)-X_-(h_N(t))).
\end{equation*}
Let us call $\Delta x(t):=x^*(t)-X_-(h_N(t))\ge 0$, then, being
$X_-(h_N(t))$ a not decreasing function for $0\le t \le \frac \pi
2N^{\frac 2 3}$, there exists  $c>0$ such that
\begin{equation*}
\frac d {dt} \, \Delta x(t)\le -c\, \Delta x(t)-\frac
{d}{dt}\,X_-(h_N(t))\le -c\,\Delta x(t), \quad \Delta x(\T_{N,T}'')=\d
\end{equation*}
hence $\Delta x(t)\le \d e^{-c  (t-\T_{N,T}'')}$  for any
$\T_{N,T}'' \le t \le \frac \pi 2 N^{2/3}$, then follows the
result. \qed

\vskip.6cm

{\bf Proof of Proposition \ref{prop:escape}.} The proof consists
of three steps. \vskip.1cm

\emph{Step I.} We prove, at first, that there exists $c>0$ such
that, for any $\ga>0$ small enough,
\begin{equation}\label{eq:B2(2)}
\p^-_{N,T}\left\{|m(\T_{N,T}')-x^*(\T_{N,T}')|\le  N^{-\frac 1
3+\frac \gamma 2} \right\}\ge 1-cN^{-\gamma}
\end{equation}

\vskip.1cm

 We have $|x^*(t)-x_c|\le \d$ for $\T_{N,T} \le t \le
\T_{N,T}'$, then there exists $c_0>0$ independent of $\d$ such
that
\begin{eqnarray*}
0\le \frac {\partial}{\partial x}\,F(x^*(t),h_N(t)) \le \b x_c
(1+c_0 \d)(x_c-x^*(t))+cT^2 N^{-2/3}
\end{eqnarray*}

thus, in particular,  there exists $c>0$  such that,
$\p^-_{N,T}$-a.s., for any $\T_{N,T} \le t \le \T_{N,T}' \wedge
\hat\T_{N,T}$,
\begin{equation*}
|\psi(t)|=\psi(t)\le c[(x_c-x^*(t))+T^2N^{-\frac 2 3}]\le
c[(x_c-\hat x(t))+T^2N^{-\frac 2 3}],
\end{equation*}

 the last inequality descending from
\eqref{esc2}. For $C_0$ as in Lemma \ref{lemma:Q1}, referring to
$\eqref{Theta}$ for the definition of
$\Theta_{\T_{N,T},\T_{N,T}'}(t)$, there exist $c,c'>0$ such that
\begin{eqnarray*}
&&\Theta_{\T_{N,T},\T_{N,T}'}(t)\le c \,N^{-1/3 +\gamma/2}\,
\exp\left\{\int_{\T_{N,T}}^t \psi(u)du\right\}\\
&&\hskip.8cm\le c \, \frac {N^{-1/3 +\gamma/2}}{\sqrt
T}\,\exp\left\{c'\,\int_{\T_{N,T}}^t
\[(x_c-\hat x(s)) +\frac{T^2} {N^{2/3}}\]\, ds\right\}\label{eq:B2(5)}
\end{eqnarray*}

Let us define $\Upsilon_{T}:=(\nu /4\b x_cT-2\mu T)$, then, by the
definition of $\hat x(t)$ in \eqref{esc2}, the exponent in
$\eqref{eq:B2(5)}$ is bounded by
\begin{eqnarray*}
c'\, \int_{\T_{N,T}}^t \(\frac 1
{\Upsilon_{T}N^{1/3}-t}+\frac{T^2}{N^{2/3}}\)\, ds =c'\,
\frac{T^2}{N^{ 2/ 3}}\,(t-\T_{N,T}) \ln
\frac{|\Upsilon_{T}N^{1/3}-t|}{|\Upsilon_{T}N^{1/3 }-\T_{N,T}|}
\end{eqnarray*}
Since, by Lemma \ref{Q1}, $\T_{N,T}'< \T_{N,T}+C_0\, {N^{1/3}}/T$,
$\p^-_{N,T}$-a.s., there exist $c,c',c''>0$ such that
\begin{eqnarray*}
\sup_{\T_{N,T} \le t \le \T_{N,T}' \wedge \hat\T_{N,T}}
\Theta_{\T_{N,T},\T_{N,T}'}(t) \le c\,  N^{-1 /3 + \gamma/
2}\,\frac{t-\Upsilon_{T}N^{1/3}}{\T_{N,T}-\Upsilon_{T}N^{1/3}}\;
e^{c'T N^{-1/3}}\\\le c'' \,N^{-1/3 + \gamma /2}
\end{eqnarray*}
for $T$ large enough, thus, by $\ref{eq:B2(6)}$,
\begin{equation*}
\p^-_{N,T}\bigg\{\sup_{\T_{N,T} \le t \le \T_{N,T}'\wedge
\hat\T_{N,T}} |m(t)-x^*(t)|\le c \,N^{-1/3+\gamma/2} \bigg\}\ge
1-c\,N^{-\gamma}
\end{equation*}

in particular, with the same probability $\T_{N,T}'<\hat\T_{N,T}$,
thus \eqref{eq:B2(2)} follows.

\vskip.1cm

\emph{Step II.} We prove, now, that there exists $c>0$ such that,
for any $\ga>0$ small enough,
\begin{equation}\label{eq:B2(3)}
\p^-_{N,T}\left\{|m(\T_{N,T}'')-x^*(\T_{N,T}'')|\le c N^{-1/3
+\gamma/2}\right\}\ge 1-cN^{-\gamma}
\end{equation}

\vskip.1cm

 We have  $|\partial F(x^*(t),h_N(t))/\partial x|\le \max\{1,\b-1\}:=c_\b$,
 thus, by \eqref{stepII}, there exists $c>0$ such that
\begin{equation*}
 \sup_{\T_{N,T}' \le t \le \T_{N,T}'' \wedge \hat\T_{N,T}}\bigg|\int_{\T_{N,T}'}^t \psi(u)du \bigg|\le
 c
 \end{equation*}
 We can use the same arguments of Step I, there exists $c>0$ such
 that
\begin{equation*}
\sup_{\T_{N,T}'\le t \le \T_{N,T}'' \wedge
\hat\T_{N,T}}\Theta_{\T_{N,T}',\T_{N,T}''}(t)\le
c\,\(|m(\T_{N,T}')-x^*(\T_{N,T}')|+N^{-(1-\gamma)/2}\)
\end{equation*}
$\p^-_{N,T}$-a.s., thus, by $\eqref{eq:B2(6)}$ and
$\eqref{eq:B2(2)}$, we have
\begin{equation*}
\p^-_{N,T}\bigg\{\sup_{\T_{N,T}' \le t \le \T_{N,T}'' \wedge
\hat\T_{N,T}}|m(t)-x^*(t)|\le c\, N^{-1 /3 +\gamma/2}\bigg\}\ge
1-c\,N^{-\gamma}
\end{equation*}
thus \eqref{eq:B2(3)} follows since  $\T_{N,T}''<\hat\T_{N,T}$
with the same probability.

\vskip.1cm

\emph{Step III.} We conclude the proof of the Proposition.  We
have $|x^*(t)-X_-(h_N(t))|\le \d$ for $t \ge \T_{N,T}''$, thus,
for small $\d$,
\begin{equation*}
\frac {\partial}{\partial x}F(x^*(t),h_N(t))=(1+\OO(\d))\, \frac
{\partial}{\partial x}F(X_-(h_N(t)),h_N(t))
\end{equation*}
On the other hand, there exists $c>0$ such that  $\frac
{\partial}{\partial x}F(X_-(h_N(t)),h_N(t))\le -c$, for any
$\T_{N,T}'' \le t \le \frac \pi 2\, N^{2/3}$, hence there exists
$c'>0$ such that
\begin{equation*}
 \sup_{\tau_2 \le t \le \frac \pi 2 N^{\frac 2
3}\wedge \hat\tau}|\psi(t)|\le -c'
\end{equation*}

Let us fix $\T_{N,T}'' \le t_* \le \frac \pi 2 \,N^{2/ 3}\wedge
\hat\T_{N,T}$, thus
\begin{eqnarray*}
\Theta_
{\T_{N,T}'',t_*}(t)=|m(\T_{N,T}'')-x^*(\T_{N,T}'')|\,e^{\int_{\T_{N,T}''}^{t}\psi(u)du}+
4 \sqrt{t} N^{-(1-\ga)/2}
\end{eqnarray*}
then, by \eqref{eq:B2(3)},
\begin{equation}\label{eq:B2(7)}
\p^-_{N,T}\left\{\Theta_ {\T_{N,T}'',t_*}(t) \le N^{-1 /3+ \gamma/
2}e^{-
 c'(t_*-\T_{N,T}'')}+ 4 \sqrt{t_*} N^{-(1-\ga)/2}\right\}\ge 1-c
\,N^{-\gamma}
\end{equation}
Let us fix, now, $\mu''>\mu'>\mu$ and choose $t_*=\mu''T N^{1/3}$,
thus, by Lemma \ref{lemma:Q1}and Lemma \ref{lemma:Q2},
$t_*>\T_{N,T}''$, $\p^-_{N,T}$-a.s., hence, by $\eqref{eq:B2(6)}$
and $\eqref{eq:B2(7)}$ we get
\begin{equation}\label{p}
\p^-_{N,T}\bigg\{\sup_{\T_{N,T}'' \le t\le t_* \wedge \hat
\T_{N,T}}|m(t)-x^*(t)|\le T N^{-\frac 1 3+\frac \gamma
2}\bigg\}\ge 1-cN^{-\gamma}
\end{equation}
then, in particular,  with the same probability $\hat \T_{N,T}
>t_*$. We have $\T_{N,T}''< \mu'T N^{1/3}<t_*$ $\p^-_{N,T}$-a.s.,
then
\begin{equation}\label{Q30}
\p^-_{N,T}\left\{|m(\mu' T N^{1/3})-x^*(\mu' T N^{1/3})|\le  T
N^{-\frac 1 3+\frac \gamma 2}\right\}\ge 1-cN^{-\gamma}
\end{equation}
On the other hand, by Lemma \ref{lemma:Q3}, $\T_{N,T}''<\mu'T
N^{1/3}$ $\p^-_{N,T}$-a.s., hence
\begin{equation}\label{Q31}
\p^-_{N,T}\left\{|m(\mu' T N^{1/3})-x^*(\mu' T N^{1/3})|\le
N^{-\frac 1 2}\right\}=1
\end{equation}
 thus
\eqref{eq:escape} follows from \eqref{Q30} and \eqref{Q31}. \qed

\section{Behavior far from criticalities}
   \label{app:D}

In this Section we give some results concerning the dynamics in
the stable region. Theorem \ref{thm:A20} provides a law for the
behavior of $m_N(t)$ in $N^{2/3}[-\frac \pi 2,-\eta]$ and
$N^{2/3}[\eta , \frac \pi 2]$, $\eta>0$. Recall that $\p_{N}$ is
the probability law of $m_N(t)$ in $N^{2/3}[-\frac \pi 2,\frac \pi
2]$ given $m_N(-\frac \pi 2 N^{2/3})=m_N^0$. For any fixed $\eta \in
[-\frac \pi 2, \frac \pi 2]$, we denote by $\p_{N}^\eta$ the law
of $m_N(t)$ in $N^{2/3}[\eta , \frac \pi 2]$ given $m_N(\eta
N^{2/3})=m_N^0$. For $\mathcal H^\pm_\g(I)$, $I \subseteq
\mathbb{R}$, as  in \eqref{H+-}, we prove the following result.

\vskip.6cm

      \begin{thm}
         \label{thm:A20}
For any $\eta,\gamma >0$ small enough and $\ga'>\ga>0$, if $|m_N^0-X_+(0)|\le N^{-\frac 1
2+\gamma}$ then
           \begin{equation}
           \label{eq:C1}
\lim_{N \to \i}\;\p_{N}\Bigg\{\mathcal H^+_{\ga'} \( N^{\frac 2 3}\[-\frac
\pi 2,-\eta\]\)\Bigg\}=1.
          \end{equation}
For any $\eta,\gamma >0$ small enough and $\ga'>\ga>0$, if $|m_N^0-X_+(h_N(\eta N^{\frac 2 3}))|\le
N^{-\frac 12 +\ga}$ then
\begin{equation}
           \label{eq:C2}
\lim_{N \to \i}\;\p_{N}^\eta\Bigg\{ \mathcal H^+_{\ga'}\(
N^{\frac 2 3}\[\eta , \frac \pi 2\]\)\Bigg\}=1.
          \end{equation}
For  any $\eta \in [-\frac \pi 2, \frac \pi 2)$, $\gamma'>\gamma$, if $|m_N^0-X_-(h_N(\eta N^{\frac 2 3}))|\le
N^{-\frac 1 2+\ga}$ then
\begin{equation}
           \label{eq:C3}
\lim_{N \to \i}\;\p_{N}^\eta\Bigg\{ \mathcal H^-_{\ga'}\(
N^{\frac 2 3}\[\eta , \frac \pi 2\]\)\Bigg\}=1.
          \end{equation}

      \end{thm}

Theorem \ref{thm:AA19} provides  a connection between the critical
and the stable regions.

\vskip.6cm

\begin{thm}
           \label{thm:AA19}
 There is $c>0$  so that for any $T$ large enough, $\gamma,\eta,\eps>0$
         \begin{equation}\label{eq:A19'}
\limsup_{N \to \i} \p_{N}\left\{|Y_N(-T)-  T|\ge \epsilon  \:
\Big| \:  \mathcal H^+_{\ga}\(\{-\eta N^{2/3}\}\)\right\}\le e^{-c
\eps^2 T}
         \end{equation}
         and
\begin{equation}\label{eq:A19''}
\limsup_{N \to \i} \p_N\left\{\(\mathcal H^+_{\ga}\(\{\eta
N^{2/3}\}\)\)^c\: \Big|\:
  |Y_N(T)-T|\le  \eps \right\}\le e^{-c \eps^2T}
         \end{equation}

               \end{thm}

\vskip.3cm

For the proof of Theorem \ref{thm:A20} and \ref{thm:AA19} see
Section 2.5 in \cite{Ca}.

\section{Conclusion of the Proof of the main result}
\label{section:Concl}

At this stage Theorem \ref{thm:TEO} is an almost direct
consequence of Theorem \ref{thm:A19}, Proposition \ref{prop:3.5}
and Proposition \ref{prop:3.6}, that we are going to prove.

 \vskip.6cm

{\bf Proof of Theorem \ref{thm:A19}.} Let us fix $\ga,\eps>0$
small enough. Recalling  that $\p_N$ is the law of $m(t)$ with
$m(- \pi N^{\frac 23}/2)=m_N^0$, suppose $|m_N^0-X_+(h_N(0))|\le
N^{-1/2+\gamma}$, then,  for any fixed $\eta>0$,
\begin{eqnarray*}
\p_N\big\{ |Y(-T)-T|> \eps \big\}\le \p_N\left\{\(\mathcal
H^+_{2\ga}\(\{-\eta N^{2/3}\}\)\)^c\right\} +\;\p_N\left\{
|Y(-T)-T|> \eps \: \Big| \: \mathcal H^+_{2\ga}\(\{-\eta
N^{2/3}\}\)\right\}
\end{eqnarray*}
thus the result follows from \eqref{eq:TEO1} and \eqref{eq:A19'}.
\qed

\vskip.6cm

{\bf Proof of Proposition \ref{prop:3.5}.} For $\p^{*}_{N,-T,y}$
and $\p^*_{-T,y}$ as defined in Section \ref{section:LD},
$\p_{N,-T,y}\{Y_N \in \e^\pm_{T}\}=\p^*_{N,-T,y}\{Y_N \in
\e^\pm_{T}\}$ and $\p_{-T,y}\{Y \in \e^\pm_{T}\}=\p^*_{-T,y}\{Y
\in \e^\pm_{T}\}$, thus Proposition \ref{prop:3.5} follows
directly from Proposition \ref{prop:z,x}. \qed

\vskip.6cm

{\bf Proof of Proposition \ref{prop:3.6}.} For any
$\eta,\gamma>0$,  $\mathcal H_{\ga}^\pm(I)$, $I \subseteq
\mathbb{R}$, as in \eqref{H+-}, we have
\begin{eqnarray*}
&&\hskip-.5cm \p_N\Big\{ \(\mathcal H^{+}_{\ga}\(N^{2/3}\[\eta,
\frac \pi
2\]\)\)^c\: \big| \: Y_N \in \e^+_{T} \Big\} \\
&&\le \p_N\left\{ \(\mathcal H^{+}_{\ga}\(N^{2/3}\[\eta, \frac \pi
2\]\)\)^c\: \Big| \: \mathcal H^+_{ \ga /2}\(\{\eta
N^{2/3}\}\)\right\}\\
 && \hskip1cm +\,\p_N\left\{  \(\mathcal H^+_{\ga/ 2}\(\{\eta
N^{2/3}\}\)\)^c\: \Big| \:   |Y_N(T)-T| \le \eps \right\}
\end{eqnarray*}
then the plus case of \eqref{m} follows from \eqref{eq:A19''} and
\eqref{eq:C2}. Analogously, for any $\eta, \ga, \mu'>\mu$
independent of $N$,   we have
\begin{eqnarray*}
&&\hskip-.5cm\p_N\left\{ \(\mathcal H^{-}_{\ga}\(N^{2/3}\[\eta,
\frac \pi 2\]\)\)^c\: \big| \: Y_N \in\e^-_{T} \right\}
\\
&& \le \p_N\left\{\(\mathcal H^{-}_{\ga}\(N^{2/3}\[\eta, \frac \pi
2\]\)\)^c\:\Big| \:   \mathcal H^-_{ \ga /2}\(\{\mu' T N^{1/3}\}\)\right\}\\
 && \hskip1cm+\,\p_N\left\{ \(\mathcal H^-_{\ga/ 2}\(\{\mu' T N^{1/3}\}\)\)^c\: \Big| \:   Y_N \in\e^-_{T} \right\}
\end{eqnarray*}
thus the minus case of \eqref{m} follows from \eqref{eq:escape}
and \eqref{eq:C3}, since  $\mu' T N^{-1/3}<<\eta$ for large $N$.
\qed

\vskip.6cm

\begin{lemma}\label{lemma:T}
We have
\begin{equation}\label{T5}
\lim_{T \to \i}\lim_{N \to \i}\Big|\p_N\big\{Y_N
\in\e^\pm_{T}\big\}-\mathbf 1_{|y-T|\le
\eps}\;\p_{N,-T,y}\big\{Y_N \in \e^\pm_{T}\big\}\Big|=0
\end{equation}
and
\begin{equation}\label{T5'}
\lim_{T \to \i}\lim_{N \to \i}\Big|\p_N\big\{Y_N \in\e^+_{T}\cup
\e^-_T\big\}-\mathbf 1_{|y-T|\le \eps}\;\p_{N,-T,y}\big\{Y_N \in
\e^+_{T}\cup \e^-_T\big\}\Big|=0.
\end{equation}

\end{lemma}

 {\bf Proof.} We prove only \eqref{T5}.We show at first that, for any fixed
$y:|y-T|\le
 \eps$, $\eps>0$  small enough,
\begin{equation}\label{Q15}
\lim_{T \to \i}\lim_{N \to \i}\Big|\p_N\left\{ Y_N \in
\e^\pm_{T}\: \big| \: |Y_N(-T)-T|\le \eps \right\}-
\p_{N,-T,y}\big\{ Y_N \in\e^\pm_{T}\big\}\Big|=0.
\end{equation}

We have
\begin{eqnarray*}
\inf_{|y-T|\le \eps} \p_{N,-T,y}\big\{ \e^\pm_{T}\big\}\le
\p_N\left\{  \e^\pm_{T}\: \big| \: |Y(-T)-T|\le \eps \right\}\le
\sup_{|y-T|\le \eps}\p_{N,-T,y}\big\{ \e^\pm_{T}\big\}
\end{eqnarray*}
thus, in order to prove \eqref{Q15}, it is sufficient to show
that, for any couple $y,\bar y:|y-T|,|\bar y-T|\le \eps$,
\begin{equation}\label{T3}
\lim_{T \to \i}\lim_{N \to \i}\big|\p_{N,-T,y}\big\{
\e^\pm_{T}\big\}-\p_{N,-T,\bar y}\big\{ \e^\pm_{T}\big\}\big|=0
\end{equation}
\eqref{T3} follows since, for  $Y(t),\bar Y(t)$ solutions of
\eqref{x} starting at $-T$ respectively from $y, \bar y$, by
Proposition \ref{prop:3.6} we have
\begin{equation}\label{Q8}
\lim_{N \to \i}\big|\p_{N,-T,y}\{ Y_N \in \e^\pm_{T}\}-\p_{-T,y}\{
Y \in\e^\pm_{T}\}\big|=0
\end{equation}
and,  by Proposition \ref{lemma:Y1},
\begin{equation}\label{Q9}
\lim_{T\to \i}\big|\p_{-T,y}\{ Y \in \e^\pm_{T}\}-\p_{-T,\bar y}\{
\bar Y \in\e^\pm_{T}\}\big|=0
\end{equation}
thus \eqref{T3} follows from \eqref{Q8} and \eqref{Q9}. We have,
now
\begin{eqnarray*}
\Big|\p_N\big\{Y_N \in \e^\pm_{T}\big\}-\p_N\big\{Y_N \in
\e^\pm_{T}\: \big| \: |Y_N(-T)-T|\le\eps \big\}\Big| \le 2 \,
\p_N\big\{|Y_N(-T)-T|>\eps\big\}. \label{Q7}
\end{eqnarray*}

From Theorem \ref{thm:A19} we know that the term in \eqref{Q7} is
vanishingly small for large $T$,  then \eqref{T5} directly follows
from \eqref{Q15}. \qed

\vskip.6cm

{\bf Conclusion of the proof of Theorem \ref{thm:TEO}.} We just
need to prove \eqref{eq:TEO2} since
 the proof of \eqref{eq:TEO1} has been proved in Section  \ref{app:D} as a part of Theorem \ref{thm:A20} (see
 \eqref{eq:C1}).\\
Let us suppose $|m_N^0-X_+(0)|\le N^{1/2+\ga}$. We have
\begin{eqnarray}
\Big|\p_N\left\{ \mathcal H_{\ga}^\pm\(N^{2/3}\[\eta, \frac \pi
2\]\)\right\}-\p_N\left\{Y_N \in
\e^\pm_{T}\right\}\Big|\nonumber\\\le \p_N\left\{\( \mathcal
H_{\ga}^\pm\(N^{2/3}\[\eta, \frac \pi 2\]\)\)^c\: \Big| \:
Y_N \in \e^\pm_{T}\right\}\nonumber\\
 +
 \p_N\left\{\( \mathcal H_{\ga}^\mp\(N^{2/3}\[\eta, \frac \pi
2\]\)\)^c \: \Big| \: Y_N \in\e^\mp_{T}\right\}\label{Q10} \\+
 \p_N\big\{ Y_N \notin \e_{T}^+ \cup \e_{T}^- \big\}\label{Q11}
\end{eqnarray}
From \eqref{eq:3.16'}, \eqref{eq:xx??} and \eqref{T5'} we have
\begin{eqnarray*}
\lim_{T \to \i}\lim_{N \to \i} \p_N\big\{ Y_N \notin \e_{T}^+ \cup
\e_{T}^- \big\}=0.
\end{eqnarray*}

From Proposition \ref{prop:3.6} we know that the terms in
\eqref{Q10} and \eqref{Q11} are vanishingly small for large $T$
and $N$, thus, from \eqref{T5} we have
\begin{equation}\label{Q13}
\lim_{T\to \i}\lim_{N \to \i}\Big|\p_N\left\{ \mathcal
H_{\ga}^\pm\(N^{2/3}\[\eta, \frac \pi 2\]\)\right\}-\mathbf
1_{|y-T|\le \eps}\;\p_{N,-T,y}\big\{Y_N \in
\e^\pm_{T}\big\}\Big|=0.
\end{equation}
Suppose $|y-T|\le \eps$, $Y(t)$ as in Proposition \ref{prop:3.6},
then
\begin{eqnarray}\label{Q14}
\hskip-.5cm \Big|\p_{N,-T,y}\big\{Y_N \in \e^\pm_{T}\big\}-p_\pm
\Big| \le \Big| \p_{N,-T,y}\big\{Y_N \in \e^\pm_{T} \big\}-
\p_{-T,y}\big\{Y \in\e^\pm_{T}\big\} \Big|\\
  +\;
\Big| \p_{-T,y}\big\{Y \in \e^\pm_{T}\big\}-p_\pm \Big|\nonumber
\end{eqnarray}
then \eqref{eq:TEO2} follows from \eqref{Q13}, \eqref{Q14},
\eqref{eq:3.15} and  \eqref{eq:3.16}. \qed

\appendix
\section{Appendix}

In  this paper we mainly make use of techniques of comparison with
Gaussian Processes. In this Appendix we provide some Gaussian
Inequalities and a comparison Lemma.

\vskip.6cm

\paragraph{Marcus-Shepp inequality for Gaussian processes.}
There is a classical result of Landau and Shepp \cite{LS} and
Marcus and Shepp \cite{MS} that gives an estimate on the
probability for a  general centered Gaussian process of escaping
from a large ball. If $G(t)$ is an a.s. bounded, centered Gaussian
process of variance $\sigma^2(t)$, then
\begin{equation}\label{MS'}
\lim_{\la \to \i}\frac 1 {\la^{2}}\ln \p \left\{\sup_{t \in I}
G(t) \ge \la\right\}=-\frac 1 {2 \sigma^2_I} \quad \text{with}
\quad \sigma^2_I:=\sup_{t \in I}\sigma^2(t)
\end{equation}
An almost immediate consequence of \eqref{MS'} is that for any
$\la$ large enough, $\d$ small enough,
\begin{equation}\label{MS}
\p \left\{\sup_t \frac {|G(t)|}{\sigma(t)} \ge \la\right\}\le 2
e^{-\frac {\la^2}2(1-\d)}.
\end{equation}

\vskip.6cm

\paragraph{Small Deviations for Gaussian Markov Processes.}
We give a result of Li (see \cite{Li}) dealing with  the
probability, for a Gaussian Markov process, of  escaping from a
small ball. Let $G(t)$ be a continuous centered  Gaussian Markov
process of covariance $\sigma(s,t)\neq 0$ for $t_0<s<t<t_1$. We
can write $\sigma(s,t)=G(s)H(t)$ with $G,H>0$ and $G/H$ non
decreasing on $(t_0,t_1)$, then
\begin{equation}\label{SB'}
\lim_{\eps\to 0} \eps^2 \ln \p\left\{  \sup_{t_0< t \le
t_1}|G(t)|<\eps \right\}=-\frac {\pi^2}8
\int_{t_0}^{t_1}(G'H-H'G)dt.
\end{equation}

We apply \eqref{SB'} to processes of the kind
\begin{equation}
G(t)=\int_{t_0}^t e^{-\int_u^t a(s)\,ds}\;dw_u, \quad t_0\le t \le
t_1,
\end{equation}
 we get
\begin{equation}\label{SB}
\lim_{\eps\to 0} \eps^2 \log \p\left\{  \sup_{t_0< t \le
t_1}|G(t)|<\eps \right\}=-\frac {\pi^2}8
\(1-e^{-\int_{t_0}^{t_1}a(s)\,ds}\).
\end{equation}

\vskip.6cm

\paragraph{Comparison with Gaussian Processes.}

In the thesis we repeatedly make use of a comparison argument
comparing the solution of a linear SDE with the solution of a more
general SDE, let us see.

Let $G_t$ be a solution of the problem
\begin{equation}\label{XXX}
 dG_t=(a(t)G_t+ b(t))dt+\xi dw_t,
\end{equation}
with $ a, b: \mathbb{R}^+ \ra \mathbb{R}$ bounded on bounded
intervals and $\xi \in \mathbb{R}$, then $G(t)$ is a Gaussian
process of the form
\begin{equation*}
G(t)=G(t_0)\,e^{\int_{t_0}^t a(s)\, ds}+\int_{t_0}^t b(s)\,
e^{\int_s^t a(u)\, du}\, ds + \xi \int_{t_0}^t e^{\int_s^t
a(u)\,du}dw_s.
\end{equation*}

Consider, now, the processes $v(t)$ solution of
\begin{equation*}
 dv_t= c(v_t,t)dt +\xi dw_t
\end{equation*}
with the same noise of \eqref{XXX},   $c: \mathbb{R} \times
\mathbb{R}^+ \ra \mathbb{R}$ globally Lipschitz.

\vskip.6cm

\begin{lemma}\label{lemma:comp}
For $G(t),v(t)$ as above we define
$\d_t:=c(G_t,t)-[a(t)G_t+b(t)]$, $\Delta_t:=G_t-v_t$, and let
$\tau \in \mathbb{R}^+$ be a generic random variable. Suppose
\begin{equation*}
 \sgn (\Delta_\tau) = \sgn
(\d_\tau) \quad \text{ or }\quad \Delta_\tau=0,
\end{equation*}
then
\begin{equation*}
\sgn (\Delta_t) = \sgn (\d_t) \quad \text{for any} \quad  \tau \le
t \le \inf\{s \ge \tau: \d_s=0\} \quad \text{a.s.}
\end{equation*}
\end{lemma}

{\bf Proof.} We have
\begin{equation*}
d\Delta_t=(a(t)\Delta_t+\d_t)dt
\end{equation*}
thus, for any $\tau \ge 0$
\begin{equation*}
\Delta(t)=\Delta(\tau)\, e^{\int_\tau^t  a(s)\,
ds}+\int_\tau^t\d(s)\, e^{\int_s^t a(u)\, du}\, ds
\end{equation*}
then follows the result. \qed

\vskip1cm

{\bf Acknowledgments.} I am grateful to my advisor Errico Presutti
for having suggested the problem  and for relevant discussions. I
 aknowledge kind hospitality at the Department of Mathematics,
University of Roma Tor Vergata. I finally thank the Department of Mathematics of the
  University of Modena and Reggio Emilia for partial financial support.

\end{document}